\documentclass[hidelinks,onefignum,onetabnum]{siamart251216}

\usepackage{lipsum}
\usepackage{amsfonts}
\usepackage{amsmath}
\usepackage{amssymb}
\usepackage{mathtools}

\usepackage{graphicx}
\usepackage{epstopdf}

\usepackage{algorithm}
\usepackage{algpseudocode}

\usepackage{empheq}
\usepackage{enumitem}

\usepackage{graphicx}
\usepackage{subcaption}

\usepackage{xcolor}
\usepackage{multirow}

\hypersetup{
	colorlinks=true,
	citecolor=blue,
	linkcolor=blue,
	urlcolor=blue
}

\ifpdf
\DeclareGraphicsExtensions{.eps,.pdf,.png,.jpg}
\else
\DeclareGraphicsExtensions{.eps}
\fi

\newsiamremark{remark}{Remark}

\newsiamremark{hypothesis}{Hypothesis}
\crefname{hypothesis}{Hypothesis}{Hypotheses}

\newsiamthm{claim}{Claim}

\newsiamremark{fact}{Fact}
\crefname{fact}{Fact}{Facts}

\newsiamremark{assumption}{Assumption}
\crefname{assumption}{Assumption}{Assumptions}
\Crefname{assumption}{Assumption}{Assumptions}

\def\RR{\mathbf{R}}

\headers{Legendre--Galerkin Method for Kinetic Opinion Dynamics}{R. Wang, Y. Xiong}

\title{A Legendre--Galerkin Method for Kinetic Models of Opinion Dynamics with Quasi-Invariant Consistency\thanks{Submitted to the editors DATE.
		\funding{This research was supported by the National Natural Science Foundation of China (Nos. 12571413). Y. X. was partially supported by Tang Scholar at Beijing Normal University.}}}

\author{Rui Wang\thanks{School of Mathematical Sciences, Beijing Normal University, Beijing 100875, China 
		(\email{rwang0913@mail.bnu.edu.cn}).}
	\and Yunfeng Xiong\thanks{School of Mathematical Sciences and Laboratory of Mathematics and Complex Systems, Ministry of Education, Beijing Normal University, Beijing 100875, China (\email{yfxiong@bnu.edu.cn}) (The corresponding author).}}

\usepackage{amsopn}

\ifpdf
\hypersetup{
  pdftitle={A Legendre--Galerkin Method for Kinetic Models of Opinion Dynamics with Quasi-Invariant Consistency},
  pdfauthor={Rui Wang and Yunfeng Xiong}
}
\fi

\begin{document}

\maketitle

\begin{abstract}
Kinetic models of opinion dynamics provide a mesoscopic framework for  describing collective behavior arising from interactions
among a large number of individuals.
Highly accurate numerical approximation to such Boltzmann-type equations, however,
remains challenging due to the nonlinear inelastic collision
operators involving multiple integrals coupling opinion
variables and random perturbations. 
To address the difficulties, we develop a Legendre--Galerkin spectral
method for a general class of kinetic models with opinion-dependent
interaction rules.
Our method exploits an exact decoupling
of the weak collision kernel to obtain a time-independent
tensor representation of the Galerkin collision operator.
A recursive construction reduces collision-tensor assembly to
one-dimensional quadrature, allowing the coefficients
to be precomputed for  subsequent use.
A rigorous $L^2$-error estimate, along with mass preservation,  is established for the fully discrete scheme under suitable
regularity and stability assumptions.
In the quasi-invariant regime, we derive the limiting
Fokker--Planck equation and establish consistency between
the exact Legendre--Galerkin operators of the rescaled
Boltzmann equation and its limit for each fixed polynomial degree.
Numerical experiments support the theoretical results,
illustrate quasi-invariant consistency, and capture
the formation of opinion clusters under the bounded-confidence interactions.
\end{abstract}



\section{Introduction}
\label{sec:introduction}

Understanding how collective behavior emerges from interactions
among individuals is a fundamental problem in the mathematical
study of social systems
\cite{weidlich2003sociodynamics,castellano2009statistical}.
A central question in opinion dynamics is how social interactions
and information exchange shape individual opinions
\cite{acemoglu2011opinion}.
Kinetic models provide a mesoscopic framework for studying
the evolution of opinion distributions in large populations
\cite{albi2017recent,toscani2006kinetic}, explaining how different interaction mechanisms can lead to consensus,
polarization or the persistence of distinct opinion groups
\cite{hegselmann2002opinion,bennaim2003bifurcations,motsch2014heterophilious,wang2017noisy},
and consequently have attracted considerable attention. Toscani~\cite{toscani2006kinetic} introduced a Boltzmann-type
formulation combining interpersonal compromise with random
self-thinking effects.
A related kinetic model was developed and mathematically
analyzed by Boudin and Salvarani~\cite{boudin2009kinetic}.
Subsequent extensions account for agent heterogeneity through
interactions between ordinary agents and strong opinion
leaders~\cite{during2009leaders}, or through personal conviction
that modulates compromise and self-thinking
\cite{brugna2015conviction}.
The kinetic framework also accommodates richer opinion spaces
and interaction structures.
Boudin et al.~\cite{boudin2010multidimensional} formulated
a multidimensional model for choices among multiple alternatives,
while D\"uring and Wolfram~\cite{during2015inhomogeneous}
incorporated leadership and spatial dependence into
inhomogeneous models.
On evolving networks, Albi et al.~\cite{albi2017network}
coupled binary opinion exchange with the creation and removal
of connections, and developed numerical methods to capture
the long-time behavior of the resulting system.
The transition from opinions to preferences was addressed by
Pareschi et al.~\cite{pareschi2019hydrodynamic}, where
hydrodynamic models of opinion-driven preference formation were derived.

In spite of the apparent differences in their interaction rules, the above models share very similar mechanisms at the microscopic level, namely, the binary opinion exchanges with
opinion-dependent compromise weights and noise amplitudes.
For two agents with opinions $w,w_{\ast}\in I=[-1,1]$,
the general interaction rule takes the form
\begin{equation}
	\begin{aligned}
		w'
		&=\alpha(w)w+\beta(w)w_{\ast}+ \eta d(w) ,\\
		w_{\ast}'
		&=\alpha(w_{\ast})w_{\ast}+\beta(w_{\ast})w + \eta_{\ast} d(w_{\ast}).
	\end{aligned}
	\label{eq:general_binary_interaction}
\end{equation}
Here $w',w_{\ast}'$ are the corresponding post-interaction opinions.
The functions $\alpha,\beta:I\to[0,1]$ are the self-retention
and compromise weights and satisfy
$\alpha(w)+\beta(w)=1$ for every $w\in I$.
Thus, in the absence of random fluctuations, each
post-interaction opinion is a convex combination of the
two pre-interaction opinions.
The function $d:I\to[0,+\infty)$ determines the amplitude
of the random fluctuations.
The random variables $\eta$ and $\eta_{\ast}$ are independent
with common probability law $\mu$ and support
$B\subset\mathbb R$.
This gives rise to a general class of opinion kinetic equations.
Let $f(w,t)\geq0$ denote the opinion density on $I$,
with initial datum $f_0\geq0$ satisfying
$\int_I f_0(w)\,\mathrm dw=1$.
For $\varphi\in C(I)$, define
$\langle f,\varphi\rangle
:=\int_I f(w,t)\varphi(w)\,\mathrm dw$.
The Boltzmann equation in weak form reads
\begin{equation}
	\frac{\mathrm{d}}{\mathrm{d}t}\langle f,\varphi\rangle
	=
	\bigl\langle Q(f,f),\varphi\bigr\rangle
	=
	\nu\int_I\int_I
	\underbrace{\mathbb{E}_\mu
		\bigl[\varphi(w')-\varphi(w)\bigr]}_{\textup{collision kernel}}
	f(w,t)f(w_{\ast},t)
	\,\mathrm{d}w_{\ast}\,\mathrm{d}w,
	\label{eq:general_boltzmann_equation}
\end{equation}
where $\nu>0$ is the interaction frequency and
$\mathbb E_\mu$ denotes the expectation under a probability measure $\mu$.
Taking $\varphi\equiv1$ in
\eqref{eq:general_boltzmann_equation} shows that
the total mass is conserved.

Solving \eqref{eq:general_boltzmann_equation} numerically
is challenging because the nonlinear inelastic collision
operator involves multiple integrals over interacting
opinions and averaging over random perturbations.
The opinion-dependent interaction rule couples these
variables within the collision integral
\cite{toscani2006kinetic}.
Direct simulation Monte Carlo (DSMC) methods can approximate the kinetic
equations by simulating the underlying binary interactions
\cite{pareschi2001montecarlo,pareschi2013multiagent} and have
been applied to leader--follower models
\cite{during2009leaders} and network-dependent dynamics
\cite{albi2017network}.
But their accuracy is affected by statistical fluctuations
\cite{caflisch1998monte}, and an accurate reconstruction
of the opinion density may require a large number of
particles \cite{albi2017network}.
Another approach uses Fokker--Planck approximations
derived from kinetic opinion models in the
quasi-invariant interaction regime
\cite{toscani2006kinetic,boudin2009quasi,during2009leaders,albi2017network}.
Although structure-preserving discretizations are available
for Fokker--Planck equations
\cite{chang1970practical,pareschi2018structure,duan2021structure},
this approximation relies on a suitable scaling of weak
interactions and small random perturbations.

In order to obtain highly accurate numerical solutions, the usage of spectral methods seems to be a natural choice, yet remains less explored for kinetic models of opinion dynamics due to the complexity in the  collision operator.
Indeed, spectral methods have been successfully applied
to the classical Boltzmann equation. The Fourier--Galerkin methods 
provide spectrally accurate
approximations of sufficiently smooth solutions of the
truncated, periodized problem
\cite{pareschi2000boltzmann},
with stability and convergence established under
suitable assumptions
\cite{filbet2011analysis,hu2021stability}.
Fast spectral methods exploit the collision geometry and
separable approximations of the Fourier weights to express
the discrete collision operator as sums of convolutions
that can be evaluated efficiently using FFTs
\cite{mouhot2006fast,filbet2006nlogn,gamba2017fast}. 
Nonetheless,  their extension to
opinion-dependent interactions is not trivial as opinion-dependent compromise weights and noise amplitudes
need not preserve the convolution structure required
by these fast algorithms.
Moreover, the bounded opinion interval is not naturally
periodic.
A direct Fourier basis expansion may therefore lose
spectral accuracy if the periodic extension of the
opinion density lacks sufficient regularity.
Evaluating the gain term in strong form presents another
difficulty, as it may require reconstructing the
pre-interaction opinions and computing the corresponding
Jacobian factor \cite{toscani2006kinetic}.
The weak formulation avoids it by
evaluating test functions directly at post-interaction
opinions.
These observations motivate us to adopt the Galerkin approach
\cite{gamba2018galerkin,kitzler2019polynomial,wang2019hermite,hu2022petrov}, and Legendre polynomials provide a natural basis for this
formulation, as they are orthogonal with respect to the
unit weight on $I$ and require no periodic extension
of the opinion density.

In this work, we develop a Legendre--Galerkin spectral
method for \eqref{eq:general_boltzmann_equation}.
The key to resolve the complexity of the collision operator is to exploit an exact finite decoupling of the weak collision
kernel, combined with Legendre orthogonality, to obtain a concise
tensor decomposition.
We construct the tensor coefficients recursively
using noise moments and the Legendre three-term recurrence to avoid explicit expansion of high-degree composite
polynomials.
This construction reduces tensor assembly to one-dimensional
Gauss--Legendre quadrature.
Since the coefficients are time-independent, they can be
precomputed, thereby facilitating subsequent collision evaluations.
A rigorous $L^2$-error estimate is established for the fully discrete scheme under the stated
regularity and stability assumptions.
The consistency between
the exact Legendre--Galerkin operators of the rescaled
Boltzmann equation and its Fokker--Planck limit in the quasi-invariant regime is also proved for each fixed polynomial degree.
Numerical experiments support our theoretical results,
illustrate quasi-invariant consistency, and capture
the formation of opinion clusters under the bounded-confidence interactions.

The rest of the paper is organized as follows.
Section~\ref{sec:kinetic_model} relates the general
interaction rule to representative models and examines
preservation of the opinion interval and the
quasi-invariant limit.
Section~\ref{sec:legendre_method} develops the fully
discrete Legendre--Galerkin spectral method and its
extension to bounded-confidence interactions.
Section~\ref{sec:numerical_analysis} analyzes errors
and establishes quasi-invariant consistency
for the exact Galerkin operators.
Section~\ref{sec:numerical_experiments} reports
numerical experiments, and Section~\ref{sec:conclusion}
concludes the paper.

\section{Kinetic formulation of opinion dynamics}
\label{sec:kinetic_model}

Building on the kinetic formulation in
Section~\ref{sec:introduction}, we relate
\eqref{eq:general_binary_interaction}
to representative models from the literatures and derive
the Fokker--Planck equation in the quasi-invariant regime.

\subsection{Representative models}
\label{subsec:binary_interaction_rule}

The binary interaction \eqref{eq:general_binary_interaction}
is \emph{admissible} if $w',w_{\ast}'\in I$ whenever
$w,w_{\ast}\in I$ and $\eta,\eta_{\ast}\in B$.
Suppose that $B\subset[-\eta_0,\eta_0]$ for some
$\eta_0\geq0$.
Using $\alpha(w)+\beta(w)=1$, we obtain
$|w'|\leq
\alpha(w)|w|+\beta(w)+\eta_0d(w)
=
1-\alpha(w)(1-|w|)+\eta_0d(w)$.
The same estimate applies to $w_{\ast}'$ with $w$ replaced
by $w_{\ast}$.
Thus, a sufficient condition for admissibility is
\begin{equation}
	\eta_0d(w)\leq\alpha(w)(1-|w|),
	\qquad w\in I.
	\label{eq:general_admissibility_condition}
\end{equation}

The following examples illustrate
\eqref{eq:general_binary_interaction}
and its bounded-confidence extension.

\begin{enumerate}[
	label=(\roman*),
	leftmargin=2.4em,
	labelsep=0.6em,
	topsep=0.4\baselineskip,
	itemsep=0.3\baselineskip,
	parsep=0pt,
	partopsep=0pt
	]
	
	\item
	The opinion update introduced by
	Toscani~\cite{toscani2006kinetic} reads
	\begin{equation*}
		w'
		=
		\bigl(1-\gamma P(|w|)\bigr)w
		+
		\gamma P(|w|)w_{\ast}
		+
		\eta D(|w|).
	\end{equation*}
	This corresponds to $\alpha(w)=1-\gamma P(|w|)$,
	$\beta(w)=\gamma P(|w|)$, and $d(w)=D(|w|)$.
	Here $\gamma$ is the compromise strength, while $P$ and $D$
	describe the opinion dependence of compromise and random
	fluctuations, respectively.
	
	\item
	Denoting the attraction function by $\rho$, the interaction
	rule of Boudin and Salvarani~\cite{boudin2009kinetic}
	can be written as
	\begin{equation*}
		w'
		=
		\frac{1+\rho(w)}{2}w
		+
		\frac{1-\rho(w)}{2}w_{\ast}.
	\end{equation*}
	The corresponding coefficients are
	$\alpha(w)=(1+\rho(w))/2$,
	$\beta(w)=(1-\rho(w))/2$, and $d(w)=0$.
	Larger values of $\rho(w)$ correspond to greater retention
	of the initial opinion and weaker attraction toward the
	pairwise average.
	
	\item
	In the single-class specialization of the model of
	Monica and Bergenti~\cite{monica2017analytic},
	the opinion update takes the form
	\begin{equation*}
		w'
		=
		(1-\gamma)w+\gamma w_{\ast},
	\end{equation*}
	where $\gamma$ is compromise weight, with $\alpha(w)=1-\gamma$,
	$\beta(w)=\gamma$, and $d(w)=0$.
	
	\item
	In the model of Brugna and Toscani~\cite{brugna2015conviction},
	each agent is characterized by an opinion $w\in I$ and
	a conviction $x\geq0$. The opinion update is
	\begin{equation*}
		w'
		=
		\bigl(1-\gamma\Psi(x)P(|w|)\bigr)w
		+
		\gamma\Psi(x)P(|w|)w_{\ast}
		+
		\eta\Phi(x)D(|w|).
	\end{equation*}
	For each fixed $x$, this has the form
	\eqref{eq:general_binary_interaction} with
	$\alpha(w)=1-\gamma\Psi(x)P(|w|)$,
	$\beta(w)=\gamma\Psi(x)P(|w|)$, and
	$d(w)=\Phi(x)D(|w|)$.
	The functions $\Psi$ and $\Phi$ describe how conviction
	modifies the propensity to compromise and the amplitude
	of random fluctuations, respectively.

	\item
	Compromise may occur only between agents whose opinions
differ by at most a prescribed threshold $\Delta$
\cite{albi2015kinetic,albi2017network,pareschi2019hydrodynamic},
 with the interaction rule
\begin{equation}
	w'
	=
	w+\theta(w)\chi_{\{|w-w_{\ast}|\leq\Delta\}}(w_{\ast}-w)
	+\eta d(w),
	\qquad 0<\Delta\leq2,
	\label{eq:bounded_confidence_interaction}
\end{equation}
and the corresponding rule for $w_{\ast}'$ obtained
by exchanging $w$ with $w_{\ast}$ and $\eta$ with $\eta_{\ast}$.
Here $\theta:I\to[0,1]$ is the compromise weight and $\chi$ is the indicator function.

\end{enumerate}

\subsection{Quasi-invariant limit}
\label{subsec:fp_limit}

The quasi-invariant regime of
\eqref{eq:general_binary_interaction} is analyzed following
the classical argument in~\cite{toscani2006kinetic}. For $0<\varepsilon\ll1$, we introduce
the scaling $\beta_\varepsilon(w)=\varepsilon b(w)$,
$\alpha_\varepsilon(w)=1-\varepsilon b(w)$, and
$\eta_\varepsilon=\sqrt{\lambda\varepsilon}\, \eta$, where $\lambda>0$, $b:I\to[0,+\infty)$ is independent
of $\varepsilon$, and $\eta$ is a (normalized) bounded random variable
with probability law $\mu$ satisfying
$\mathbb E_{\mu}[\eta]=0$ and $\mathbb E_{\mu}[\eta^2]=1$. The corresponding interaction rule is
\begin{equation}
	\begin{aligned}
		w_\varepsilon'
		&=
		w+\varepsilon b(w)(w_{\ast}-w)
		+\sqrt{\lambda\varepsilon}\,\eta d(w),
		\\
		w_{\varepsilon,*}'
		&=
		w_{\ast}+\varepsilon b(w_{\ast})(w-w_{\ast})
		+\sqrt{\lambda\varepsilon}\,\eta_{\ast}d(w_{\ast}),
	\end{aligned}
	\label{eq:qi_scaled_interaction}
\end{equation}
where $\eta$ and $\eta_{\ast}$ are independent and identically distributed.

Let $f_\varepsilon(w,t)$ be a solution
of the Boltzmann equation associated with
\eqref{eq:qi_scaled_interaction}, with
$f_\varepsilon(\cdot,0)=f_0$.
With the macroscopic time $\tau=\nu\varepsilon t$
and $g_\varepsilon(w,\tau)
=f_\varepsilon(w,\frac{\tau}{\nu\varepsilon})$, the following proposition identifies the limiting equation.

\begin{proposition}[Quasi-invariant limit]
	\label{prop:quasi_invariant_limit}
	Suppose the following conditions hold:
	\begin{itemize}
	\item[(1)] $b,d\in C(I)$ with $b,d\geq0$;
	
	\item[(2)] There exists
$\varepsilon_0\in(0,1)$ such that
$\varepsilon_0\|b\|_{L^\infty(I)}\leq1$ and
\begin{equation}
	\sqrt{\lambda\varepsilon}\,
	\|\eta \|_{L^\infty}d(w)
	\leq
	\bigl(1-\varepsilon b(w)\bigr)(1-|w|),
	\qquad
	w\in I,\quad 0<\varepsilon\leq\varepsilon_0.
	\label{eq:qi_admissibility}
\end{equation}
\end{itemize}
Then $g_\varepsilon$ admits a weakly convergent subsequence to a family of probability
	measures $g$ on $I$ as $\varepsilon\to0$.
	Every such limit \(g\) satisfies the weak formulation
	\eqref{eq:fp_weak_form} of the Fokker--Planck equation
	\begin{equation}
		\partial_\tau g
		=
		-\partial_w
		\left[
		b(w)\bigl(m(\tau)-w\bigr)g
		\right]
		+
		\frac{\lambda}{2}
		\partial_w^2
		\left[d(w)^2g\right],
		\qquad w\in(-1,1),
		\label{eq:fp_limit}
	\end{equation}
	where $m(\tau)=\int_I w\,g(\tau,\mathrm dw)$,
	with initial datum
	$g(0,\mathrm dw)=f_0(w)\,\mathrm dw$.
\end{proposition}

\begin{remark}
The nonnegativity of $b$ and the bound
$\varepsilon_0\|b\|_{L^\infty(I)}\leq1$
ensure that $0\leq\alpha_\varepsilon(w),\beta_\varepsilon(w)\leq1$
for all $w\in I$ and $0<\varepsilon\leq\varepsilon_0$.
Condition \eqref{eq:qi_admissibility} guarantees
admissibility by \eqref{eq:general_admissibility_condition}.
\end{remark}

\begin{proof}
	The rescaled weak formulation reads
\begin{equation}
	\frac{\mathrm d}{\mathrm d\tau}
	\int_I\varphi(w)g_\varepsilon(w,\tau)\,\mathrm dw
	=
	\frac{1}{\varepsilon}
	\int_{I^2}
	\mathbb{E}_{\mu}\left[
	\varphi(w_\varepsilon')-\varphi(w)
	\right]
	g_\varepsilon(w,\tau)g_\varepsilon(w_{\ast},\tau)
	\,\mathrm dw\,\mathrm dw_{\ast}.
	\label{eq:qi_rescaled_weak_boltzmann}
\end{equation}
Set
$\delta_\varepsilon:=w_\varepsilon'-w
=\varepsilon b(w)(w_{\ast}-w)+\sqrt{\lambda\varepsilon}\,\eta d(w)$.
Since \eqref{eq:qi_scaled_interaction} is admissible and $I$ is
convex, the line segment joining $w$ and $w_\varepsilon'$ is contained
in $I$. Thus, for $\varphi\in C^3(I)$, Taylor's formula gives
\begin{equation*}
	\varphi(w_\varepsilon')-\varphi(w)
	=
	\varphi'(w)\delta_\varepsilon
	+\frac12\varphi''(w)\delta_\varepsilon^2
	+r_{\varphi,\varepsilon},
	\quad
	|r_{\varphi,\varepsilon}|
	\leq
	\frac16\|\varphi'''\|_{L^\infty(I)}
	|\delta_\varepsilon|^3.
\end{equation*}
Using $\mathbb E_{\mu}[\eta]=0$ and
$\mathbb E_{\mu}[\eta^2]=1$, we obtain
\begin{equation*}
	\frac{1}{\varepsilon}\mathbb{E}_{\mu}[\delta_\varepsilon]
	=
	b(w)(w_{\ast}-w),
	\quad
	\frac{1}{\varepsilon}\mathbb{E}_{\mu}[\delta_\varepsilon^2]
	=
	\lambda d(w)^2
	+
	\varepsilon b(w)^2(w_{\ast}-w)^2.
\end{equation*}
Moreover, since $I$ is bounded, $b,d\in C(I)$, and $\eta$ is bounded, we have
$\mathbb{E}_{\mu}[|\delta_\varepsilon|^3]\leq C\varepsilon^{3/2}$
uniformly in $(w,w_{\ast})\in I^2$. Consequently,
\begin{equation}
	\frac{1}{\varepsilon}
	\mathbb{E}_{\mu}\left[
	\varphi(w_\varepsilon')-\varphi(w)
	\right]
	=
	b(w)(w_{\ast}-w)\varphi'(w)
	+\frac{\lambda}{2}d(w)^2\varphi''(w)
	+R_{\varphi,\varepsilon}(w,w_{\ast}),
	\label{eq:fp_kernel_expansion}
\end{equation}
where the Taylor remainder $R_{\varphi,\varepsilon}$  satisfies $\|R_{\varphi,\varepsilon}\|_{L^\infty(I^2)}
\leq
C_\varphi\sqrt{\varepsilon},
~~
0<\varepsilon\leq\varepsilon_0$.  

Substituting \eqref{eq:fp_kernel_expansion} into
\eqref{eq:qi_rescaled_weak_boltzmann} and using the unit mass of
$g_\varepsilon$ gives
\begin{equation}
	\begin{aligned}
		\frac{\mathrm d}{\mathrm d\tau}
		\int_I\varphi(w)g_\varepsilon(w,\tau)\,\mathrm dw
		&=
		\int_I
		b(w)\bigl(m_\varepsilon(\tau)-w\bigr)
		\varphi'(w)g_\varepsilon(w,\tau)\,\mathrm dw
		\\
		&\quad+
		\frac{\lambda}{2}
		\int_I
		d(w)^2\varphi''(w)g_\varepsilon(w,\tau)\,\mathrm dw
		+\mathcal R_{\varphi,\varepsilon}[g_{\varepsilon}(\tau),g_{\varepsilon}(\tau)],
	\end{aligned}
	\label{eq:fp_weak_expansion_gamma}
\end{equation}
where $m_\varepsilon(\tau):=\int_I w g_\varepsilon(w,\tau)\,\mathrm dw$, and
\begin{equation*}
	\mathcal R_{\varphi,\varepsilon}[g_{\varepsilon}(\tau),g_{\varepsilon}(\tau)]
	:=
	\int_{I^2}
	R_{\varphi,\varepsilon}(w,w_{\ast})
	g_\varepsilon(w,\tau)g_\varepsilon(w_{\ast},\tau)
	\,\mathrm dw\,\mathrm dw_{\ast}.
\end{equation*}
The unit mass of $g_\varepsilon$ implies
$|\mathcal R_{\varphi,\varepsilon}[g_{\varepsilon}(\tau),g_{\varepsilon}(\tau)]|
\leq C_\varphi\sqrt{\varepsilon},~0<\varepsilon\leq\varepsilon_0$.
	
	Since $|m_\varepsilon|\leq1$ and $b,d$ are bounded,
\eqref{eq:fp_weak_expansion_gamma} implies
$\left|
\frac{\mathrm d}{\mathrm d\tau}
\int_I\varphi(w)g_\varepsilon(w,\tau)\,\mathrm dw
\right|
\leq C_\varphi$
uniformly in $\varepsilon$.
Consequently, for every $T>0$,
\[
\left|
\int_I\varphi(w)
\bigl(g_\varepsilon(w,t)-g_\varepsilon(w,s)\bigr)
\,\mathrm dw
\right|
\leq C_\varphi|t-s|,
\qquad s,t\in[0,T].
\]
Moreover, the unit mass of $g_\varepsilon$ gives
$\left|
\int_I\varphi(w)g_\varepsilon(w,t)\,\mathrm dw
\right|
\leq \|\varphi\|_{L^\infty(I)}$, $t\in[0,T]$.
	Applying the Arzel\`a--Ascoli theorem to a countable
	dense subset of $C(I)$ contained in $C^3(I)$,
	together with compactness of probability measures on $I$
	and a diagonal argument, we obtain a subsequence
	$\varepsilon_n\to0$ such that $g_{\varepsilon_n}$
	converges to a family of probability measures $g$.
	In particular, $m_{\varepsilon_n}\to m$ uniformly
	on such intervals, where
	$m(\tau):=\int_I w\,g(\tau,\mathrm dw)$.
	
	By integrating \eqref{eq:fp_weak_expansion_gamma}
	over $[0,\tau]$ and passing to the limit and using the continuity of $b,d$ and the remainder estimate, it gives
	\begin{equation}
		\begin{aligned}
			\int_I\varphi(w)\,g(\tau,\mathrm dw)
			&-\int_I\varphi(w)f_0(w)\,\mathrm dw=
				\frac{\lambda}{2}
			\int_0^\tau\int_I
			d(w)^2\varphi''(w)
			\,g(s,\mathrm dw)\,\mathrm ds
			\\
			&\qquad+\int_0^\tau\int_I
			b(w)\bigl(m(s)-w\bigr)\varphi'(w)
			\,g(s,\mathrm dw)\,\mathrm ds.
		\end{aligned}
		\label{eq:fp_weak_form}
	\end{equation}
	This identity holds for every $\varphi\in C^3(I)$.
	Taking $\varphi\in C_c^\infty((-1,1))$ shows that
	\eqref{eq:fp_limit} holds in the sense of distributions
	in the interior of $I$.
\end{proof}

\section{Legendre-Galerkin spectral method}
\label{sec:legendre_method}

This section presents a Legendre--Galerkin spectral method
for \eqref{eq:general_boltzmann_equation}.
The exact decoupling of the weak collision kernel produces
a time-independent collision tensor.
Its coefficients can be precomputed using a recursive
construction and one-dimensional Gauss--Legendre quadrature.
An extension to bounded-confidence interactions will also be discussed.

Throughout this paper, we assume that $\mu$ has finite absolute moments up to
order $N$, and set $M_q:=\mathbb E_\mu[\eta^q]$,
$0\leq q\leq N$.
These moments depend only on the prescribed noise law
and are precomputed analytically or by numerical
quadrature.

\subsection{Legendre--Galerkin semi-discretization}
\label{subsec:legendre_galerkin}

Let $L_k$ denote the Legendre polynomial of degree $k$ on
$I$, with $L_0\equiv1$ and $L_1(w)=w$. The polynomials satisfy
\begin{equation}
	\int_I
	L_k(w)L_\ell(w)\,\mathrm dw
	=
	h_k\delta_{k\ell},
	\qquad
	h_k:=\frac{2}{2k+1}.
	\label{eq:legendre_orthogonality}
\end{equation}
For a truncation degree \(N\geq1\), set
$V_N:=\operatorname{span}
\{L_0,L_1,\ldots,L_N\}$
and seek an approximate opinion density of the form
\begin{equation}
	f_N(w,t)
	=
	\sum_{k=0}^{N}
	\widehat f_{k}(t)L_k(w).
	\label{eq:legendre_expansion}
\end{equation}

The Legendre--Galerkin semi-discretization seeks
$f_N(\cdot,t)\in V_N$ such that
\begin{equation}
	\frac{\mathrm{d}}{\mathrm{d}t}
	\langle f_N,\varphi_N\rangle
	=
	\bigl\langle Q(f_N,f_N),\varphi_N\bigr\rangle,
	\qquad
	\forall\,\varphi_N\in V_N.
	\label{eq:galerkin_weak_form}
\end{equation}
Taking $\varphi_N=L_\ell$ and using
\eqref{eq:legendre_orthogonality}, we obtain $	h_\ell
\frac{\mathrm{d}\widehat f_\ell(t)}{\mathrm{d}t}
=
\bigl\langle Q(f_N,f_N),L_\ell\bigr\rangle,
~
0\leq\ell\leq N$.
It remains to evaluate the projected collision terms
$\bigl\langle Q(f_N,f_N),L_\ell\bigr\rangle$.

\subsection{Finite decoupling of the weak collision kernel}
\label{subsec:kernel_decoupling}

Under the interaction rule \eqref{eq:general_binary_interaction},
the variable $w_{\ast}$ enters $w'$ only through the term
$\beta(w)w_{\ast}$ linearly depending on $w_{\ast}$, which allows an exact finite expansion of the
weak collision kernel in $w_{\ast}$.

\begin{lemma}[Finite expansion in $w_{\ast}$]
	\label{lem:finite_kernel_expansion}
	Let $0\leq\ell\leq N$ be fixed. Define
	\begin{equation}
		\mathcal K_\ell(w,w_{\ast})
		:=
		\mathbb E_\mu\!\left[L_\ell(w')-L_\ell(w)\right],\quad 0\leq\ell\leq N.
		\label{eq:weak_kernel}
	\end{equation}
	Then, for every fixed $w\in I$, the function
	$w_{\ast}\mapsto\mathcal K_\ell(w,w_{\ast})$ is a polynomial of degree at
	most $\ell$. Consequently, there exist unique coefficients
	$\Phi_{\ell m}(w)$, $0\leq m\leq\ell$, such that
	\begin{align}
		\mathcal K_\ell(w,w_{\ast})
		&=
		\sum_{m=0}^{\ell}
		\Phi_{\ell m}(w)L_m(w_{\ast}),
		\label{eq:kernel_expansion} \\
		\Phi_{\ell m}(w)
		&=
		\frac{1}{h_m}
		\int_I
		\mathcal K_\ell(w,z)L_m(z)\,\mathrm dz,
		\qquad
		0\leq m\leq\ell.
		\label{eq:Phi_definition}
	\end{align}
	In particular, all Legendre coefficients of
	$\mathcal K_\ell(w,\cdot)$ of order greater than $\ell$ vanish.
\end{lemma}

\begin{proof}
	Fix $w\in I$. Since $L_\ell$ is a polynomial of degree $\ell$,
	write
	$L_\ell(x)=\sum_{r=0}^{\ell}c_{\ell r}x^r$. Using
	$w'=\alpha(w)w+\beta(w)w_{\ast}+\eta d(w)$, the multinomial formula
	gives
	\begin{equation*}
		\begin{aligned}
			\mathbb{E}_\mu\left[L_\ell(w')\right]
			&=
			\sum_{r=0}^{\ell}c_{\ell r}
			\sum_{\substack{i,j,q\geq0\\i+j+q=r}}
			\frac{r!}{i!j!q!}
			\alpha(w)^i\beta(w)^j d(w)^q
			M_q\,w^i w_{\ast}^j.
		\end{aligned}
	\end{equation*}
	The right-hand side is a polynomial in $w_{\ast}$ of degree at most
	$\ell$. Since the loss term $L_\ell(w)$ is independent of
	$w_{\ast}$, the same conclusion holds for
	$\mathcal K_\ell(w,w_{\ast})$.
	
	Because $\{L_0,\ldots,L_\ell\}$ is a basis of
	$\mathbb P_\ell(I)$, the expansion
	\eqref{eq:kernel_expansion} is unique. Multiplying it
	by $L_n(w_{\ast})$ and integrating over $I$, we obtain
	\eqref{eq:Phi_definition} by the orthogonality relation
	\eqref{eq:legendre_orthogonality}. The vanishing of all coefficients of
	order greater than $\ell$ follows directly from the
	orthogonality of $L_n$ to $\mathbb P_\ell(I)$ for $n>\ell$.
\end{proof}

Substituting \eqref{eq:kernel_expansion} into
\eqref{eq:general_boltzmann_equation} with $f=f_N$ and
$\varphi=L_\ell$, and using \eqref{eq:legendre_expansion}
and \eqref{eq:legendre_orthogonality}, we obtain
\begin{equation}
	\begin{aligned}
		\bigl\langle Q(f_N,f_N),L_\ell\bigr\rangle
		&=
		\nu
		\sum_{m=0}^{\ell}
		\left(
		\int_I
		L_m(w_{\ast})f_N(w_{\ast},t)\,\mathrm{d}w_{\ast}
		\int_I
		\Phi_{\ell m}(w)f_N(w,t)\,\mathrm{d}w
		\right)
		\\
		&=
		\nu
		\sum_{m=0}^{\ell}\sum_{k=0}^{N}
		A_{\ell mk}
		h_m\widehat f_m(t)
		\widehat f_k(t),
		\quad 0\leq\ell\leq N,
	\end{aligned}
	\label{eq:Q_separated}
\end{equation}
where the time-independent collision tensor is defined by
\begin{equation}
	A_{\ell mk}
	:=
	\int_I
	\Phi_{\ell m}(w)L_k(w)\,\mathrm{d}w,
	\qquad
	0\leq m\leq\ell\leq N,
	\quad
	0\leq k\leq N.
	\label{eq:A_definition}
\end{equation}
Combining with \eqref{eq:Q_separated},
we obtain the Legendre--Galerkin semi-discrete system
\begin{equation}
	\begin{aligned}
		\frac{\mathrm d\widehat f_\ell(t)}{\mathrm dt}
		&=
		\frac{\nu}{h_\ell}
		\sum_{m=0}^{\ell}\sum_{k=0}^{N}
		h_m A_{\ell mk}\widehat f_m(t)\widehat f_k(t)=:R_\ell\bigl(\widehat{\mathbf f}(t), \widehat{\mathbf f}(t)\bigr),
		\qquad 0\leq\ell\leq N,
	\end{aligned}
	\label{eq:modal_ode_final}
\end{equation}
with
$\widehat{\mathbf f}(t)
:=(\widehat f_0(t),\ldots,\widehat f_N(t))^{\mathrm T}$. Equivalently, the system takes the vector form
$\frac{\mathrm d\widehat{\mathbf f}(t)}{\mathrm dt}
=\RR\bigl(\widehat{\mathbf f}(t), \widehat{\mathbf f}(t)\bigr)$, where
$\RR(\widehat{\mathbf f}, \widehat{\mathbf f}):=
\bigl(R_0(\widehat{\mathbf f}, \widehat{\mathbf f}),\ldots,
R_N(\widehat{\mathbf f}, \widehat{\mathbf f})\bigr)^{\mathrm T}$.

%
%


The collision tensor $A_{\ell mk}$ depends only on the interaction
weights $\alpha$ and $\beta$, the noise amplitude $d$, and the
moments of the noise distribution. It can therefore be assembled
once in a preprocessing stage and reused throughout the time
integration.

\subsection{Recursive computation of the collision tensor
	$A_{\ell mk}$}
\label{subsec:recursive_tensor}

For $0\leq\ell\leq N$ and $0\leq s\leq N-\ell$,
the same degree argument as in
Lemma~\ref{lem:finite_kernel_expansion} shows that,
for each fixed $w\in I$,
$w_{\ast}\mapsto\mathbb E_\mu[\eta^sL_\ell(w')]$
is a polynomial of degree at most $\ell$.
Hence there exist unique coefficients
$\Psi_{\ell m}^{(s)}(w)$, $0\leq m\leq\ell$, such that
\begin{equation}
	\mathbb E_\mu[\eta^sL_\ell(w')]
	=
	\sum_{m=0}^{\ell}
	\Psi_{\ell m}^{(s)}(w)L_m(w_{\ast}).
	\label{eq:Psi_expansion}
\end{equation}
For each admissible pair $(\ell,s)$, we set
$\Psi_{\ell m}^{(s)}=0$ for $m<0$ or $m>\ell$.

\begin{theorem}[Recursive computation of the collision tensor]
	\label{thm:A_recurrence}
	For $0\leq m\leq\ell\leq N$ and $0\leq k\leq N$,
	the collision tensor satisfies
	\begin{equation}
		A_{\ell mk}
		=
		\int_I
		\Psi_{\ell m}^{(0)}(w)L_k(w)\,\mathrm{d}w
		-
		\delta_{m0}h_\ell\delta_{\ell k}.
		\label{eq:A_recursive_formula}
	\end{equation}
	The auxiliary coefficients $\Psi_{\ell m}^{(s)}$ are initialized
	by
	\begin{subequations}
		\label{eq:Psi_initialization}
		\begin{align}
			\Psi_{00}^{(s)}(w)
			&=
			M_s,
			\quad 0\leq s\leq N,
			\label{eq:Psi_init_l0}
			\\
			\Psi_{10}^{(s)}(w)
			&=
			M_s\alpha(w)w+M_{s+1}d(w),
			\quad
			\Psi_{11}^{(s)}(w)
			=
			M_s\beta(w),
			\quad 0\leq s\leq N-1.
			\label{eq:Psi_init_l1}
		\end{align}
	\end{subequations}
The coefficients satisfy the following matrix recurrence relation
for $1\leq\ell\leq N-1$ and $0\leq s\leq N-\ell-1$:
	\begin{equation}
		\begin{aligned}
			\begin{pmatrix}
				\Psi_{\ell+1,0}^{(s)}\\
				\Psi_{\ell+1,1}^{(s)}\\
				\vdots\\
				\Psi_{\ell+1,N}^{(s)}
			\end{pmatrix}
			={}&
			\frac{2\ell+1}{\ell+1}
			\begin{pmatrix}
				\alpha(w)w
				& \dfrac{\beta(w)}{3}
				& 0
				& \cdots
				& 0
				\\[1ex]
				\beta(w)
				& \alpha(w)w
				& \dfrac{2\beta(w)}{5}
				& \ddots
				& \vdots
				\\[1ex]
				0
				& \dfrac{2\beta(w)}{3}
				& \alpha(w)w
				& \ddots
				& 0
				\\
				\vdots
				& \ddots
				& \ddots
				& \ddots
				& \dfrac{N\beta(w)}{2N+1}
				\\[1ex]
				0
				& \cdots
				& 0
				& \dfrac{N\beta(w)}{2N-1}
				& \alpha(w)w
			\end{pmatrix}
			\begin{pmatrix}
				\Psi_{\ell,0}^{(s)}\\
				\Psi_{\ell,1}^{(s)}\\
				\vdots\\
				\Psi_{\ell,N}^{(s)}
			\end{pmatrix}
			\\[1ex]
			&-
			\frac{\ell}{\ell+1}
			\begin{pmatrix}
				\Psi_{\ell-1,0}^{(s)}\\
				\Psi_{\ell-1,1}^{(s)}\\
				\vdots\\
				\Psi_{\ell-1,N}^{(s)}
			\end{pmatrix}
			+
			\frac{2\ell+1}{\ell+1}d(w)
			\begin{pmatrix}
				\Psi_{\ell,0}^{(s+1)}\\
				\Psi_{\ell,1}^{(s+1)}\\
				\vdots\\
				\Psi_{\ell,N}^{(s+1)}
			\end{pmatrix}.
		\label{eq:Psi_recurrence}
		\end{aligned}
	\end{equation}
\end{theorem}


\begin{proof}
	Since $L_0\equiv1$, relation \eqref{eq:Psi_expansion} gives
	$\mathbb{E}_\mu\!\left[\eta^sL_0(w')\right]
	=\mathbb{E}_\mu\!\left[\eta^s\right]
	=M_s$, which proves
	\eqref{eq:Psi_init_l0}. Moreover, using $L_1(x)=x$ and
	\eqref{eq:general_binary_interaction}, we obtain
	\begin{equation*}
		\mathbb{E}_\mu\!\left[\eta^sL_1(w')\right]
		=M_s\alpha(w)w
		+M_s\beta(w)w_{\ast}
		+M_{s+1}d(w).
	\end{equation*}
	Comparison of the coefficients of $L_0(w_{\ast})=1$ and
	$L_1(w_{\ast})=w_{\ast}$ gives \eqref{eq:Psi_init_l1}.
	
	For $1\leq\ell\leq N-1$, the Legendre recurrence relation reads
	\begin{equation*}
		L_{\ell+1}(w')
		=
		\frac{2\ell+1}{\ell+1}w' L_\ell(w')
		-
		\frac{\ell}{\ell+1}L_{\ell-1}(w').
	\end{equation*}
Then we obtain
	\begin{equation}
		\begin{aligned}
			\mathbb{E}_\mu\!\left[\eta^sL_{\ell+1}(w')\right]
			&=
			\frac{2\ell+1}{\ell+1}
			\alpha(w)w
			\mathbb{E}_\mu\!\left[\eta^sL_\ell(w')\right]+\frac{2\ell+1}{\ell+1}
			\beta(w)w_{\ast}
			\mathbb{E}_\mu\!\left[\eta^sL_\ell(w')\right]\\
			&\quad
			+\frac{2\ell+1}{\ell+1}
			d(w)
			\mathbb{E}_\mu\!\left[
			\eta^{s+1}L_\ell(w')
			\right]-
			\frac{\ell}{\ell+1}
			\mathbb{E}_\mu\!\left[\eta^sL_{\ell-1}(w')\right].
		\end{aligned}
		\label{eq:Psi_recurrence_aux}
	\end{equation}
Using \eqref{eq:Psi_expansion},
			$
			w_{\ast}L_j(w_{\ast})
			=
			\frac{j+1}{2j+1}L_{j+1}(w_{\ast})
			+
			\frac{j}{2j+1}L_{j-1}(w_{\ast}), j\geq1$, and	$w_{\ast}L_0(w_{\ast})=L_1(w_{\ast})$,
			the coefficient of $L_m(w_{\ast})$ in
			$w_{\ast}\sum_{j=0}^{\ell}
			\Psi_{\ell j}^{(s)}(w)L_j(w_{\ast})$ is
			$\frac{m}{2m-1}\Psi_{\ell,m-1}^{(s)}(w)
			+\frac{m+1}{2m+3}\Psi_{\ell,m+1}^{(s)}(w)$.
	Comparing the coefficients of $L_m(w_{\ast})$ in
	\eqref{eq:Psi_recurrence_aux} gives
	\eqref{eq:Psi_recurrence}.
	
	Finally, setting $s=0$ in \eqref{eq:Psi_expansion} gives
	$\mathbb{E}_\mu\!\left[L_\ell(w')\right]
	=\sum_{m=0}^{\ell}
	\Psi_{\ell m}^{(0)}(w)L_m(w_{\ast})$.
	Recalling \eqref{eq:weak_kernel}, we therefore have
	\begin{equation*}
		\mathcal K_\ell(w,w_{\ast})
		=
		\sum_{m=0}^{\ell}
		\left[
		\Psi_{\ell m}^{(0)}(w)
		-
		\delta_{m0}L_\ell(w)
		\right]
		L_m(w_{\ast}).
	\end{equation*}
	Comparing with \eqref{eq:kernel_expansion}, it gives  $\Phi_{\ell m}(w)=\Psi_{\ell m}^{(0)}(w)-\delta_{m0}L_\ell(w)$.
	Substituting it into \eqref{eq:A_definition}, it proves \eqref{eq:A_recursive_formula} by the orthogonal relation  \eqref{eq:legendre_orthogonality}.
\end{proof}

The integral in \eqref{eq:A_recursive_formula} is approximated
using an $N_q$-point Gauss--Legendre quadrature rule on $I$,
with nodes $w_i$ and weights $\omega_i$.
The quadrature tensor is defined by
\begin{equation}
	A_{\ell mk}^{(N_q)}
	:=
	\sum_{i=1}^{N_q}
	\omega_i
	\Psi_{\ell m}^{(0)}(w_i)L_k(w_i)
	-
	\delta_{m0}h_\ell\delta_{\ell k},
	\quad
	0\leq m\leq\ell\leq N,
	\quad
	0\leq k\leq N.
	\label{eq:A_quadrature}
\end{equation}
The same quadrature nodes and weights are used for all tensor entries.
At each node $w_i$, relations
\eqref{eq:Psi_init_l0}--\eqref{eq:Psi_recurrence} generate all
admissible values of $\Psi_{\ell m}^{(s)}(w_i)$, while the values
$L_k(w_i)$ are computed once and reused.

Since $\Psi_{00}^{(0)}=M_0=1$, we have
$\Phi_{00}=0$ and hence $A_{00k}=0$ and $A_{00k}^{(N_q)}=0$
for $0\leq k\leq N$.
When every $\Psi_{\ell m}^{(0)}$ is a polynomial in $w$, let
$\deg_\Psi$ denote the largest polynomial degree among all
$\Psi_{\ell m}^{(0)}$, $0\leq m\leq\ell\leq N$. Since $L_k$ has
degree at most $N$, the integrand in
\eqref{eq:A_recursive_formula} has degree at most
$\deg_\Psi+N$. Consequently,
\begin{equation}
	A_{\ell mk}^{(N_q)}
	=
	A_{\ell mk}
	\quad
	\text{for all admissible }\ell,m,k,
	\quad
	\text{whenever }
	2N_q-1\geq\deg_\Psi+N.
	\label{eq:A_quadrature_exactness}
\end{equation}

\subsection{Extension to bounded-confidence interactions}
\label{subsec:bounded_confidence_method}

Now consider the interaction rule
\eqref{eq:bounded_confidence_interaction}.
Let $Q_\Delta$ denote the corresponding collision
operator in the weak formulation
\eqref{eq:general_boltzmann_equation}.
We assume that $\theta$ and $d$ are smooth,
that the interaction is admissible, and that $\mu$
has finite absolute moments up to order $N$.

For fixed $w\in I$, define
$I_\Delta(w):=[p_\Delta(w),q_\Delta(w)]$, where
$p_\Delta(w):=\max\{-1,w-\Delta\}$ and
$q_\Delta(w):=\min\{1,w+\Delta\}$.
Then
\[
w'
=
\begin{cases}
	(1-\theta(w))w+\theta(w)w_{\ast}+\eta d(w),
	& w_{\ast}\in I_\Delta(w),\\[1mm]
	w+\eta d(w),
	& w_{\ast}\notin I_\Delta(w).
\end{cases}
\]
On the confidence branch, let
$\Psi_{\ell j}^{(s)}$ denote the coefficients
generated by \eqref{eq:Psi_initialization} and
\eqref{eq:Psi_recurrence} with
$\alpha(w)=1-\theta(w)$ and $\beta(w)=\theta(w)$.
On the complementary branch, set
$H_\ell(w):=\mathbb E_\mu[L_\ell(w+\eta d(w))]$.
This function is computed using the same construction
with $\alpha(w)=1$ and $\beta(w)=0$, taking the
zeroth Legendre coefficient at $s=0$.
Thus,
\[
\mathbb E_\mu[L_\ell(w')]
=
\begin{cases}
	\displaystyle
	\sum_{j=0}^{\ell}
	\Psi_{\ell j}^{(0)}(w)L_j(w_{\ast}),
	& w_{\ast}\in I_\Delta(w),\\[2mm]
	H_\ell(w),
	& w_{\ast}\notin I_\Delta(w),
\end{cases}
\qquad 0\leq\ell\leq N.
\]
We can introduce the restricted
Legendre inner products
\[
G_{jm}^{\Delta}(w)
:=
\int_{p_\Delta(w)}^{q_\Delta(w)}
L_j(v)L_m(v)\,\mathrm dv,
\qquad 0\leq j,m\leq N.
\]
Legendre orthogonality gives
$\int_{I\setminus I_\Delta(w)}L_m(v)\,\mathrm dv
=h_0\delta_{m0}-G_{0m}^{\Delta}(w)$.

Now the time-independent collision tensor can be calculated by
\begin{equation}
	\begin{aligned}
		A_{\ell mk}^{\Delta}
		:={}&
		\frac{1}{h_m}
		\int_I L_k(w)
		\sum_{j=0}^{\ell}
		\Psi_{\ell j}^{(0)}(w)
		G_{jm}^{\Delta}(w)\,\mathrm dw
		\\
		&+	\frac{1}{h_m}
		\int_I L_k(w)
		H_\ell(w)
		\bigl(h_0\delta_{m0}-G_{0m}^{\Delta}(w)\bigr)
		\,\mathrm dw
		-\delta_{m0}h_\ell\delta_{\ell k},
	\end{aligned}
	\label{eq:bounded_confidence_tensor}
\end{equation}
for $0\leq\ell,m,k\leq N$.
The two integral terms are the gain contributions
from the confidence and complementary regions, respectively.
The final term is the loss contribution, evaluated exactly
using Legendre orthogonality.

Substituting the Legendre expansion of $f_N$ into
the weak collision integral gives
\begin{equation}
	\bigl\langle Q_\Delta(f_N,f_N),L_\ell\bigr\rangle
	=
	\nu\sum_{m=0}^{N}h_m\widehat f_m
	\sum_{k=0}^{N}
	A_{\ell mk}^{\Delta}\widehat f_k,
	\qquad 0\leq\ell\leq N.
	\label{eq:bounded_confidence_modal_operator}
\end{equation}
The coefficient equations are
$h_\ell\,\mathrm d\widehat f_\ell/\mathrm dt
=\langle Q_\Delta(f_N,f_N),L_\ell\rangle$,
$0\leq\ell\leq N$.
The sum over $m$ extends to $N$ because the restricted
inner products couple different Legendre degrees.

The entries $G_{jm}^{\Delta}(w)$ are computed by
analytically integrating the polynomial $L_jL_m$.
The integrals in \eqref{eq:bounded_confidence_tensor}
are approximated by an $N_q$-point Gauss--Legendre rule
on each nondegenerate subinterval of $I$ obtained
by splitting at $-1+\Delta$ and $1-\Delta$.

\subsection{Fully discrete scheme}
\label{subsec:fully_discrete_scheme}

Define $\RR^{(N_q)}$ by replacing $A_{\ell mk}$ in
\eqref{eq:modal_ode_final} with its quadrature approximation
$A_{\ell mk}^{(N_q)}$ from \eqref{eq:A_quadrature}.
The corresponding coefficient vector
$\widehat{\mathbf f}^{(N_q)}(t)$ satisfies
\begin{equation}
	\frac{\mathrm d\widehat{\mathbf f}^{(N_q)}(t)}
	{\mathrm dt}
	=
	\RR^{(N_q)}
	\bigl(\widehat{\mathbf f}^{(N_q)}(t), \widehat{\mathbf f}^{(N_q)}(t)\bigr).
	\label{eq:quadrature_modal_system}
\end{equation}
For bounded-confidence interactions, use
$A_{\ell mk}^{\Delta,(N_q)}$
and extend the sum over $m$ to $N$.
Let $t_n=n\Delta t$,
$N_t=\lfloor T/\Delta t\rfloor$, and let $\widehat{\mathbf f}^{\,n}$ approximate $\widehat{\mathbf f}^{(N_q)}(t_n)$.
The classical fourth-order Runge--Kutta (RK4) method gives
\begin{equation}
	\begin{aligned}
		\mathbf k_1^{\,n}
		&=
		\RR^{(N_q)}
		\bigl(\widehat{\mathbf f}^{\,n},\widehat{\mathbf f}^{\,n}\bigr),\quad
		\mathbf k_2^{\,n}
		=
		\RR^{(N_q)}
		\left(
		\widehat{\mathbf f}^{\,n}
		+\frac{\Delta t}{2}\mathbf k_1^{\,n},
		\widehat{\mathbf f}^{\,n}
		+\frac{\Delta t}{2}\mathbf k_1^{\,n}
		\right),
		\\
		\mathbf k_3^{\,n}
		&=
		\RR^{(N_q)}
		\left(
		\widehat{\mathbf f}^{\,n}
		+\frac{\Delta t}{2}\mathbf k_2^{\,n},
			\widehat{\mathbf f}^{\,n}
		+\frac{\Delta t}{2}\mathbf k_2^{\,n}
		\right),~
		\mathbf k_4^{\,n}
		=
		\RR^{(N_q)}
		\left(
		\widehat{\mathbf f}^{\,n}
		+\Delta t\,\mathbf k_3^{\,n},
		\widehat{\mathbf f}^{\,n}
		+\Delta t\,\mathbf k_3^{\,n}
		\right),
		\\
		\widehat{\mathbf f}^{\,n+1}
		&=
		\widehat{\mathbf f}^{\,n}
		+\frac{\Delta t}{6}
		\left(
		\mathbf k_1^{\,n}
		+2\mathbf k_2^{\,n}
		+2\mathbf k_3^{\,n}
		+\mathbf k_4^{\,n}
		\right),
	\end{aligned}
	\label{eq:rk4_update}
\end{equation}
for $0\leq n<N_t$. The complete procedure is summarized in
Algorithm~\ref{alg:legendre_spectral_method}.
\begin{algorithm}[htbp]
	\caption{Legendre spectral method with recursive precomputation}
	\label{alg:legendre_spectral_method}
	\begin{algorithmic}[1]
		\Require Polynomial degree $N$; Gauss--Legendre rule
		$\{w_i,\omega_i\}_{i=1}^{N_q}$; interaction weights
		$\alpha$ and $\beta$; noise amplitude $d$;
		interaction frequency $\nu$; noise moments
		$\{M_j\}_{j=0}^{N}$; initial datum $f_0$;
		final time $T$; time step $\Delta t$.
		
		\Ensure Legendre coefficients
		$\widehat{\mathbf f}^{\,n}$, $0\leq n\leq N_t$.
		
		\Statex \textbf{Precomputation stage}
		
		\State Compute $\Psi_{\ell m}^{(s)}(w_i)$ using
		\eqref{eq:Psi_init_l0}--\eqref{eq:Psi_recurrence}.
		
		\State Assemble $A_{\ell mk}^{(N_q)}$
		using \eqref{eq:A_quadrature}.
		
		\State Set $A_{00k}^{(N_q)}=0$
		for $0\leq k\leq N$.
		
		\Statex \textbf{Time-integration stage}
		
		\State Set
		$\widehat f_k^{\,0}
		=h_k^{-1}\int_I f_0(w)L_k(w)\,\mathrm dw$
		for $0\leq k\leq N$.
		
		\For{$n=0,\ldots,N_t-1$}
		\State Compute the RK4 stages and update
		$\widehat{\mathbf f}^{\,n+1}$
		using \eqref{eq:rk4_update}.
		\EndFor
		
		\State Reconstruct $f_N^{\,n}$
		using $f_N^{\,n}(w)=\sum_{k=0}^{N}\widehat f_k^{\,n}L_k(w)$.
		
	\end{algorithmic}
\end{algorithm}

The fully discrete scheme preserves the total mass as $A_{00k}^{(N_q)}=0$ implies that $\widehat f_0^{\,n}=\widehat f_0^{\,0}$ at each RK4 stage.
By Legendre orthogonality,
\begin{equation}
	\int_I f_N^{\,n}(w)\,\mathrm dw
	=
	h_0\widehat f_0^{\,n}
	=
	h_0\widehat f_0^{\,0}
	=
	\int_I f_N^{\,0}(w)\,\mathrm dw,
	\qquad 0\leq n\leq N_t.
	\label{eq:fully_discrete_mass}
\end{equation}

\section{Numerical analysis}
\label{sec:numerical_analysis}

This section presents a fully discrete error estimate and
establishes the convergence
of the exact Legendre--Galerkin operators of the rescaled
Boltzmann equation to their Fokker--Planck counterparts
in the quasi-invariant limit.

\subsection{Error analysis}
\label{subsec:error_analysis}

The error analysis for \eqref{eq:general_binary_interaction}
covers Legendre--Galerkin approximation, collision-tensor
quadrature, and RK4 time discretization.
The noise moments are treated as exact input data throughout
the analysis.

Let $f_N$ and $\widetilde f_N$ denote the
semi-discrete approximations obtained using the exact
collision tensor $A_{\ell mk}$ and its quadrature
approximation $A_{\ell mk}^{(N_q)}$, respectively.
At each time level $t_n$, $0\leq n\leq N_t$, the total
error decomposes as
\begin{equation}
		f(\cdot,t_n)-f_N^{\,n}
		=
		\underbrace{\bigl(f(\cdot,t_n)-f_N(\cdot,t_n)\bigr)}_{\textup{Legendre
approximation}}+
		\underbrace{\bigl(f_N(\cdot,t_n)-\widetilde f_N(\cdot,t_n)\bigr)}_{\textup{collision-tensor quadrature}}
		+
		\underbrace{\bigl(\widetilde f_N(\cdot,t_n)-f_N^{\,n}\bigr)}_{\textup{RK4 time integrator}}.
	\label{eq:error_decomposition}
\end{equation}

Let $\Pi_N:L^2(I)\to V_N$ denote the $L^2(I)$-orthogonal
projection. All approximations use the same projected initial
datum, with
$f_N(\cdot,0)=\widetilde f_N(\cdot,0)=f_N^{\,0}=\Pi_Nf_0$.
To account for the numerical evaluation of the collision tensor,
we define
\begin{equation}
	\varepsilon_A
	:=
	\max_{\substack{0\leq m\leq\ell\leq N\\0\leq k\leq N}}
	\left|A_{\ell mk}-A_{\ell mk}^{(N_q)}\right|.
	\label{eq:error_tensor_quadrature}
\end{equation}

Assume that the bilinear collision operator
$Q:L^2(I)\times L^2(I)\to L^2(I)$ satisfies
the following bound, analogous to the Boltzmann
estimate in
\cite{mouhot2004regularity,filbet2011analysis,hu2021stability}:
\begin{equation}
	\|Q(u,v)\|_{L^2(I)}
	\leq C_Q\|u\|_{L^2(I)}\|v\|_{L^2(I)},
	\qquad u,v\in L^2(I),
	\label{eq:error_collision_bilinear}
\end{equation}
where $C_Q>0$ depends only on the fixed model parameters.

\begin{remark}[A sufficient condition for the bilinear estimate]	
	Estimate \eqref{eq:error_collision_bilinear} holds if
the interaction \eqref{eq:general_binary_interaction}
is admissible and $\beta(w)\geq\beta_0>0$ on $I$.
	For fixed $w\in I$ and $\eta\in B$,
	since $w'=\alpha(w)w+\eta d(w)+\beta(w)w_{\ast}$,
	a change of variables gives
	\[
	\int_I |\varphi(w')|^2\,\mathrm dw_{\ast}
	=
	\frac{1}{\beta(w)}
	\int_{\alpha(w)w+\eta d(w)-\beta(w)}^{\alpha(w)w+\eta d(w)+\beta(w)}
	|\varphi(z)|^2\,\mathrm dz
	\leq
	\beta_0^{-1}\|\varphi\|_{L^2(I)}^2
	\]
	for every $\varphi\in C(I)$.
	Applying the Cauchy--Schwarz inequality to the weak collision
	integral, we obtain
	\[
	\begin{aligned}
		\bigl|\langle Q(u,v),\varphi\rangle\bigr|
		&\leq
		\nu\Bigl(
		\beta_0^{-1/2}\|u\|_{L^1(I)}\|v\|_{L^2(I)}
		+
		\|u\|_{L^2(I)}\|v\|_{L^1(I)}
		\Bigr)\|\varphi\|_{L^2(I)}
		\\
		&\leq
		\nu\sqrt{2}(1+\beta_0^{-1/2})
		\|u\|_{L^2(I)}\|v\|_{L^2(I)}
		\|\varphi\|_{L^2(I)}.
	\end{aligned}
	\]
	The last inequality uses
	$\|z\|_{L^1(I)}\leq\sqrt{2}\|z\|_{L^2(I)}$.
	By density of $C(I)$ in $L^2(I)$ and duality,
	\eqref{eq:error_collision_bilinear} follows with
	$C_Q=\nu\sqrt{2}(1+\beta_0^{-1/2})$.
	This condition is satisfied by the binary compromise
rule, e.g., in Boudin and Salvarani
\cite[Section~2.1]{boudin2009kinetic},
	for which $\beta(w)=(1-\rho(w))/2$ and $d\equiv0$.
	Their attraction function $\rho$ is continuous on $I$
	and satisfies $0\leq\rho<1$, so that
	$\beta(w)\geq(1-\max_I\rho)/2>0$.
\end{remark}

\begin{theorem}[Fully discrete error estimate]
	\label{thm:error_fully_discrete}
Suppose that the collision operator satisfies the bilinear
estimate \eqref{eq:error_collision_bilinear}.
	Let $r\geq1$ be an integer, and let $f$ be a weak solution of
	\eqref{eq:general_boltzmann_equation} satisfying
	$f\in C^1([0,T];L^2(I))\cap C([0,T];H^r(I))$.
	Let $f_N^n$ denote the fully discrete approximation with
	$f_N^0=\Pi_N f(\cdot,0)$, and set
	$t_n=n\Delta t$ and $N_t=\lfloor T/\Delta t\rfloor$.
	Then there exist constants $N_0\geq1$,  $C_s>0$ and $C_q > 0$
	such that
	\begin{equation*}
		\begin{aligned}
			\max_{0\leq n\leq N_t}
			\|f(\cdot,t_n)-f_N^n\|_{L^2(I)}
			\leq
			C_s(N+1)^{-r}\|f\|_{C([0,T];H^r(I))}
		+C_q\varepsilon_A+C_t\Delta t^4,
		\end{aligned}
	\end{equation*}
	whenever $N\geq N_0$ and $C_q\varepsilon_A < 1$. The constant $C_s$ is independent of $N$, $N_q$, and $\Delta t$.
	The constants $C_q$ and $C_t$
 	 depend on $N$, but are independent of $N_q$ and $\Delta t$.
\end{theorem}

The proof is based on the following two lemmas.

\begin{lemma}[Legendre--Galerkin error]
	\label{thm:error_spatial}
	Let $Q$, $r$, and $f$ satisfy the assumptions of
	Theorem~\ref{thm:error_fully_discrete}.
	Then there exist an integer $N_0\geq1$ and a constant $C_s>0$
	such that, for every $N\geq N_0$, the semi-discrete solution
	$f_N$ exists uniquely on $[0,T]$ and
	\begin{equation}
		\max_{0\leq t\leq T}
		\|f(\cdot,t)-f_N(\cdot,t)\|_{L^2(I)}
		\leq
		C_s(N+1)^{-r}\|f\|_{C([0,T];H^r(I))}.
		\label{eq:error_spatial_estimate}
	\end{equation}
	The constant $C_s$ is independent of $N$.
\end{lemma}

\begin{proof}
For fixed $N$, the bilinear estimate
\eqref{eq:error_collision_bilinear} and the
$L^2(I)$-contractivity of $\Pi_N$ imply that
the right-hand side of \eqref{eq:modal_ode_final}
is locally Lipschitz in the coefficient vector.
The semi-discrete system with the prescribed
initial datum therefore admits a unique maximal
solution on $[0,T_N^{\max})$.

	Consider the projection
	error \(\rho_N:=f-\Pi_Nf\) and the discrete error 
	\(\theta_N:=\Pi_Nf-f_N\) and then
	\(f-f_N=\rho_N+\theta_N\).
	The projected initial data give \(\theta_N(\cdot,0)=0\).
	
	Subtracting \eqref{eq:galerkin_weak_form} from the
	exact weak equation tested with \(\varphi_N\in V_N\) gives
	\begin{equation}
		\langle\partial_t\theta_N,\varphi_N\rangle
		=
		\bigl\langle Q(f,f)-Q(f_N,f_N),\varphi_N\bigr\rangle,
		\qquad
		\varphi_N\in V_N.
		\label{eq:error_spatial_identity}
	\end{equation}
	
	Set $M_T:=\|f\|_{C([0,T];L^2(I))}$ and
	$T_N^*:=\min\{T,T_N^{\max}\}$, and define a cutoff time
	\begin{equation}\label{cutoff_time}
	\tau_N
	:=
	\sup\left\{
	\tau\in[0,T_N^*):
	\|\theta_N(\cdot,t)\|_{L^2(I)}\leq1
	\text{ for all }t\in[0,\tau]
	\right\}.
	\end{equation}
	Since $\theta_N(\cdot,0)=0$, continuity implies $\tau_N>0$.
	For $0\leq t<\tau_N$, the identity
	$f_N=\Pi_Nf-\theta_N$ and the contractivity of $\Pi_N$
	in $L^2(I)$ give
	$\|f_N(\cdot,t)\|_{L^2(I)}\leq M_T+1$.
    			By bi-linearity,
    			$Q(f,f)-Q(f_N,f_N)=Q(f-f_N,f)+Q(f_N,f-f_N)$.
    			Testing \eqref{eq:error_spatial_identity} with
    			$\varphi_N=\theta_N$ and applying the Cauchy--Schwarz
    			inequality and \eqref{eq:error_collision_bilinear}, we obtain
    			\[
    			\begin{aligned}
    				\frac12\frac{\mathrm d}{\mathrm dt}
    				\|\theta_N\|_{L^2(I)}^2
    				&\leq
    				C_Q
    				\bigl(\|f\|_{L^2(I)}+\|f_N\|_{L^2(I)}\bigr)
    				\|f-f_N\|_{L^2(I)}
    				\|\theta_N\|_{L^2(I)}
    				\\
    				&\leq
    				2C_Q(M_T+1)
    				\bigl(
    				\|\rho_N\|_{L^2(I)}
    				+
    				\|\theta_N\|_{L^2(I)}
    				\bigr)
    				\|\theta_N\|_{L^2(I)},
    			\end{aligned}
    			\]
	Young's inequality gives
	\begin{equation}
		\frac{\mathrm d}{\mathrm dt}
		\|\theta_N\|_{L^2(I)}^2
		\leq
		6C_Q(M_T+1)\|\theta_N\|_{L^2(I)}^2
		+
		2C_Q(M_T+1) \|\rho_N\|_{L^2(I)}^2.
		\label{eq:error_spatial_differential}
	\end{equation}
	Using \(\theta_N(\cdot,0)=0\), Gr\"onwall's inequality implies
	\begin{equation}
		\sup_{0\leq t<\tau_N}
		\|\theta_N(\cdot,t)\|_{L^2(I)}
		\leq
		C_1
		\max_{0\leq t\leq T}
		\|\rho_N(\cdot,t)\|_{L^2(I)},
		\label{eq:error_projected_spatial_estimate}
	\end{equation}
	where \(C_1>0\) depends only on \(T\), $C_Q$ and $M_T$.
	
	For every \(N\geq r-1\), applying the Legendre projection estimate
	in \cite{shen2011spectral}
	\(\alpha=\beta=0\), \(\ell=0\), and \(m=r\), we obtain
	\begin{equation}
		\max_{0\leq t\leq T}
		\|\rho_N(\cdot,t)\|_{L^2(I)}
		\leq
		C_r(N+1)^{-r}
		\|f\|_{C([0,T];H^r(I))},
		\label{eq:error_projection_estimate}
	\end{equation}
	where \(C_r>0\) is independent of \(N\). Combining
	\eqref{eq:error_projected_spatial_estimate} and
	\eqref{eq:error_projection_estimate}, we obtain
	\[
	\sup_{0\leq t<\tau_N}
	\|\theta_N(\cdot,t)\|_{L^2(I)}
	\leq
	C_1C_r(N+1)^{-r}
	\|f\|_{C([0,T];H^r(I))}.
	\]
	It suffices to choose \(N_0\geq\max\{1,r-1\}\) such that
	\(C_1C_r(N_0+1)^{-r}
	\|f\|_{C([0,T];H^r(I))}\leq1/2\). Then, for every
	\(N\geq N_0\),
	$\sup_{0\leq t<\tau_N}
	\|\theta_N(\cdot,t)\|_{L^2(I)}
	\leq\frac12$.
	
	It remains to show that $\tau_N=T$. 
	Suppose that \(\tau_N<T_N^*\). By continuity, the same bound
	holds at \(t=\tau_N\). Hence there exists \(\delta>0\) such
	that \(\tau_N+\delta<T_N^*\) and
	\(\|\theta_N(\cdot,t)\|_{L^2(I)}<1\) for
	\(t\in[\tau_N,\tau_N+\delta]\). Together with the preceding
	estimate on \([0,\tau_N]\), this contradicts the definition
	of \(\tau_N\). Therefore, \(\tau_N=T_N^*\).
	
	Moreover, we have \(T_N^{\max}>T\). Suppose that
	\(T_N^{\max}\leq T\). Since \(\tau_N=T_N^*=T_N^{\max}\), the
	preceding estimate gives
	\[
	\|f_N(\cdot,t)\|_{L^2(I)}
	\leq
	\|\Pi_Nf(\cdot,t)\|_{L^2(I)}
	+
	\|\theta_N(\cdot,t)\|_{L^2(I)}
	\leq
	M_T+1,
	\quad
	0\leq t<T_N^{\max}.
	\]
	By Legendre orthogonality and \(h_k\geq h_N>0\),
	\[
	\sum_{k=0}^{N}
	|\widehat f_k(t)|^2
	\leq
	h_N^{-1}\|f_N(\cdot,t)\|_{L^2(I)}^2
	\leq
	h_N^{-1}(M_T+1)^2.
	\]
Thus $\widehat{\mathbf f}$ remains bounded on
$[0,T_N^{\max})$.
The map $\mathbf z\mapsto\RR(\mathbf z,\mathbf z)$
is locally Lipschitz on $\mathbb R^{N+1}$.
By the ODE continuation theorem \cite{teschl2012ordinary},
the solution can therefore be extended beyond
the finite time $T_N^{\max}$.
This contradicts the maximality of $T_N^{\max}$.
Therefore, $T_N^{\max}>T$.
	It follows that \(T_N^*=T\) and \(\tau_N=T\). By continuity,
	\eqref{eq:error_projected_spatial_estimate} holds on the
	closed interval \([0,T]\). 
	
	Combining \eqref{eq:error_projected_spatial_estimate} and
	\eqref{eq:error_projection_estimate}, it arrives at \eqref{eq:error_spatial_estimate} with
	\(C_s=(1+C_1)C_r\).
\end{proof}

We next estimate the collision-tensor quadrature error.

Let $\widehat{\mathbf f}(t)$ and $\widehat{\mathbf f}^{(N_q)}(t)$
denote the coefficient vectors of $f_N(\cdot,t)$ and
$\widetilde f_N(\cdot,t)$, respectively.
For
\(\mathbf y=(y_0,\ldots,y_N)^{\mathsf T}\) and
\(\mathbf z=(z_0,\ldots,z_N)^{\mathsf T}\) in
\(\mathbb R^{N+1}\), define the weighted inner product and associated norm by
\begin{equation}
	(\mathbf y,\mathbf z)_h
	:=
	\sum_{k=0}^{N}h_ky_kz_k,
	\qquad
	\|\mathbf z\|_h^2
	:=
	\sum_{k=0}^{N}h_k|z_k|^2.
	\label{eq:error_modal_norm}
\end{equation}
By Legendre orthogonality,
\(\|\sum_{k=0}^{N}z_kL_k\|_{L^2(I)}=\|\mathbf z\|_h\).

\begin{lemma}[Quadrature error for the collision tensor]
	\label{thm:error_quadrature}
Under the assumptions of Theorem~\ref{thm:error_fully_discrete}, for each $N\geq N_0$ with $N_0$ 
 chosen in Lemma~\ref{thm:error_spatial}, there exist a constant $C_q>0$ such that the solution
$\widetilde f_N$ exists uniquely on $[0,T]$ when
$C_q \varepsilon_A < 1$, and
	\begin{equation}
		\max_{0\leq t\leq T}
		\|f_N(\cdot,t)-\widetilde f_N(\cdot,t)\|_{L^2(I)}
		\leq
		C_q\varepsilon_A,
		\label{eq:error_quadrature_estimate}
	\end{equation}
	where \(C_q\) depends on
	\(T\) and \(N\), but are independent of \(N_q\).
\end{lemma}

\begin{proof}
	Set $u:=\sum_{k=0}^{N}y_kL_k$ and
	$v:=\sum_{k=0}^{N}z_kL_k$.
	By the Legendre-Galerkin formulation, $\RR(\mathbf y, \mathbf y)$ and
	$\RR(\mathbf z, \mathbf z)$ are the coefficient vectors
	of $\Pi_NQ(u,u)$ and $\Pi_NQ(v,v)$, respectively.
	Using the identity
	$Q(u,u)-Q(v,v)=Q(u-v,u)+Q(v,u-v)$
	and the contractivity of $\Pi_N$, we obtain
	\begin{equation}
		\begin{aligned}
			\|\RR(\mathbf y, \mathbf y)
			-\RR(\mathbf z, \mathbf z)\|_h &\leq
			\|Q(u-v,u)\|_{L^2(I)}
			+\|Q(v,u-v)\|_{L^2(I)}
			\\
			&\leq
			C_Q\bigl(\|u\|_{L^2(I)}+\|v\|_{L^2(I)}\bigr)
			\|u-v\|_{L^2(I)}
			\\
			&=
			C_Q\bigl(\|\mathbf y\|_h+\|\mathbf z\|_h\bigr)
			\|\mathbf y-\mathbf z\|_h,
		\end{aligned}
		\label{eq:error_modal_quadratic_bound}
	\end{equation}
	where the second inequality uses the bilinear estimate \eqref{eq:error_collision_bilinear}, and the third equality uses  \eqref{eq:error_modal_norm}.

	Using the definition of $\varepsilon_A$ in
	\eqref{eq:error_tensor_quadrature}, the tensor
	representations of $\RR(\mathbf z,\mathbf z)$ in
	\eqref{eq:modal_ode_final} and
	$\RR^{(N_q)}(\mathbf z,\mathbf z)$ in
	\eqref{eq:quadrature_modal_system}, together with
	the Cauchy--Schwarz inequality, we obtain
	\begin{equation}
		\|\RR(\mathbf z, \mathbf z)-\RR^{(N_q)}(\mathbf z, \mathbf z)\|_h
		\leq
		D\varepsilon_A\|\mathbf z\|_h^2,
		\label{eq:error_tensor_vector_field_bound}
	\end{equation}
	where \(D>0\) depends only on \(N\) and \(\nu\).

For fixed $N_q$, the map
$\mathbf z\mapsto\RR^{(N_q)}(\mathbf z,\mathbf z)$
is quadratic and locally Lipschitz. Let $[0,T_q^{\max})$ be the maximal existence
	interval of its solution, and set
	$\mathbf e:=\widehat{\mathbf f}
	-\widehat{\mathbf f}^{(N_q)}$.
	Let
	\(M:=1+\max_{0\leq t\leq T}
	\|\widehat{\mathbf f}(t)\|_h\) and
	\(T_q^*:=\min\{T,T_q^{\max}\}\), and define a cutoff time
	\[
	\tau_q
	:=
	\sup\left\{
	\tau\in[0,T_q^*):
	\|\mathbf e(t)\|_h\leq1
	\text{ for all }t\in[0,\tau]
	\right\},
	\]
	so that \(\|\widehat{\mathbf f}^{(N_q)}(t)\|_h\leq M\) on $[0, \tau_q)$.
	Subtracting the two semi-discrete systems and taking the
	\(h\)-inner product with \(\mathbf e\), we obtain
	\begin{equation}
		\frac{\mathrm d}{\mathrm dt}\|\mathbf e\|_h^2
		\leq
		(4C_QM+1)\|\mathbf e\|_h^2
		+
		D^2M^4\varepsilon_A^2.
		\label{eq:error_quadrature_differential}
	\end{equation}
	where \eqref{eq:error_modal_quadratic_bound} and \eqref{eq:error_tensor_vector_field_bound} are used. 
	Grönwall's inequality gives
	\begin{equation}
		\sup_{0\leq t<\tau_q}
		\|\mathbf e(t)\|_h
		\leq
		C_q\varepsilon_A,
		\label{eq:error_quadrature_modal_estimate}
	\end{equation}
	where \(C_q>0\) depends on \(T\) and \(N\), but not on \(N_q\). 
	
	It remains to prove $\tau_q=T_q^*$. When $\varepsilon_A$
		is sufficiently small such that $C_q\varepsilon_A<1$,
		$\sup_{0\le t<\tau_q}\|\mathbf e(t)\|_h<1$.
		If $\tau_q<T_q^*$, continuity extends the bound
		$\|\mathbf e(t)\|_h\le1$ to $[0,\tau_q+\delta]$
		for some $\delta>0$ with $\tau_q+\delta<T_q^*$,
		which contradicts the definition of $\tau_q$.
		Hence $\tau_q=T_q^*$.
	
	Suppose \(T_q^{\max}\leq T\). Then
	$\|\widehat{\mathbf f}^{(N_q)}(t)\|_h
	\leq
	\|\widehat{\mathbf f}(t)\|_h
	+
	\|\mathbf e(t)\|_h
	\leq M$ for $0\leq t<T_q^{\max}$.
Because the map
$\mathbf z\mapsto\RR^{(N_q)}(\mathbf z,\mathbf z)$
is locally Lipschitz on $\mathbb R^{N+1}$, the solution remains bounded on $[0,T_q^{\max})$. By
the ODE continuation theorem \cite{teschl2012ordinary}, it
guarantees an extension beyond the finite time $T_q^{\max}$.
This contradicts the maximality of $T_q^{\max}$.
Therefore, $T_q^{\max}>T$.
	Thus $\tau_q=T$ and \eqref{eq:error_quadrature_modal_estimate}
	holds on $[0,T]$. 
\end{proof}

Combining the decomposition formula \eqref{eq:error_decomposition} and the above two lemmas, we complete the proof of the main theorem.
\begin{proof}[Proof of Theorem~\ref{thm:error_fully_discrete}]
	Let $N_0$ and $C_s$ be as in Lemma~\ref{thm:error_spatial},
	and fix $N\geq N_0$. For sufficiently large $N_q$ satisfying
	$C_q \varepsilon_A < 1$, Lemma~\ref{thm:error_quadrature} gives
	\[
	\max_{0\leq t\leq T}
	\|\widehat{\mathbf f}^{(N_q)}(t)\|_h
	\leq
	\max_{0\leq t\leq T}
	\|\widehat{\mathbf f}(t)\|_h
	+C_q\varepsilon_A,
	\]
	implying that the coefficient vectors
	$\widehat{\mathbf f}^{(N_q)}(t)$ are uniformly bounded
	on $[0,T]$.
	Moreover, the definition of \(\varepsilon_A\) gives
	$|A_{\ell mk}^{(N_q)}|
	\leq
	|A_{\ell mk}|+\varepsilon_A$,
	$0\leq m\leq\ell\leq N$, $0\leq k\leq N$.
	Thus $R^{(N_q)}$ and their derivatives are bounded.
	
	The standard error estimate for RK4 \cite{hairer1993solving}
	 gives $\max_{t_n\leq T}
		\left\|
		\widehat{\mathbf f}^{(N_q)}(t_n)
		-
		\widehat{\mathbf f}^{\,n}
		\right\|_h
		\leq
		C_t\Delta t^4$. By Legendre orthogonality,
	\[
	\|\widetilde f_N(\cdot,t_n)-f_N^{\,n}\|_{L^2(I)}
	=
	\left\|
	\widehat{\mathbf f}^{(N_q)}(t_n)
	-
	\widehat{\mathbf f}^{\,n}
	\right\|_h.
	\]
	Combining the decomposition \eqref{eq:error_decomposition} with
	\eqref{eq:error_spatial_estimate} and
	\eqref{eq:error_quadrature_estimate} proves
	Theorem \ref{thm:error_fully_discrete}.
\end{proof}

\subsection{Discrete quasi-invariant consistency}
\label{subsec:discrete_quasi_invariant}

Proposition~\ref{prop:quasi_invariant_limit} establishes the
quasi-invariant limit at the continuous level.
For fixed $N$, the exact Legendre--Galerkin operators
satisfy an analogous consistency property, as illustrated below.
\begin{equation*}
	\begin{array}{ccc}
		\text{Rescaled Boltzmann equation}
		&
		\xrightarrow{\ \varepsilon\to0\ }
		&
		\text{Fokker--Planck (FP) equation}
		\\
		\Big\downarrow{\scriptstyle\text{Legendre--Galerkin}}
		&&
		\Big\downarrow{\scriptstyle\text{Legendre--Galerkin}}
		\\
		\text{Approximated Boltzmann equation}
		&
		\xrightarrow[N\ \mathrm{fixed}]{\varepsilon\to0}
		&
		\text{Approximated FP equation}
	\end{array}
\end{equation*}

Let $N\geq1$ and $f_N\in V_N$ satisfy
$\int_I f_N(w)\,\mathrm dw=1$.
Set $m_N:=\int_I wf_N(w)\,\mathrm dw=h_1\widehat f_1$.
Using \eqref{eq:fp_weak_form},
we define the Fokker--Planck Galerkin terms by
\begin{equation*}
	\mathcal F_\ell(f_N)
	:=
	\int_I
	\left[
	b(w)(m_N-w)L_\ell'(w)
	+
	\frac{\lambda}{2}d(w)^2L_\ell''(w)
	\right]
	f_N(w)\,\mathrm{d}w,
	\quad 0\leq\ell\leq N.
\end{equation*}
Thus, the corresponding Legendre--Galerkin system is
$h_\ell\frac{\mathrm d\widehat f_\ell(\tau)}{\mathrm d\tau}
=\mathcal F_\ell(f_N(\cdot,\tau))$.

For the scaled interaction \eqref{eq:qi_scaled_interaction},
define the corresponding Legendre--Galerkin collision
terms by
\begin{equation*}
	\mathcal B_{\ell,\varepsilon}(f_N)
	:=
	\int_{I^2}
	\mathbb{E}_{\mu}\!\left[
	L_\ell(w_\varepsilon')-L_\ell(w)
	\right]
	f_N(w)f_N(w_{\ast})
	\,\mathrm{d}w\,\mathrm{d}w_{\ast},
	\quad 0\leq\ell\leq N.
\end{equation*}
Under the time scaling $\tau=\nu\varepsilon t$,
the corresponding Legendre--Galerkin system is
$h_\ell\frac{\mathrm d\widehat f_\ell(\tau)}{\mathrm d\tau}
=\varepsilon^{-1}\mathcal B_{\ell,\varepsilon}(f_N(\cdot,\tau)).$

\begin{theorem}[Discrete quasi-invariant consistency]
	\label{thm:qi_discrete_consistency}
	Under the same assumptions  in Proposition~\ref{prop:quasi_invariant_limit}, and let $N\geq1$ be
	fixed. Then for every $f_N\in V_N$ satisfying
	$\int_I f_N\,\mathrm dw=1$ and every $0\leq\ell\leq N$, there
	exists a constant $C_\ell>0$ independent of $\varepsilon$ and
	$f_N$, such that
	\begin{equation*}
	\left|
	\frac{1}{\varepsilon}\mathcal B_{\ell,\varepsilon}(f_N)
	-\mathcal F_\ell(f_N)
	\right|
		\leq
		C_\ell\|f_N\|_{L^1(I)}^2\sqrt{\varepsilon},
		\qquad
		0<\varepsilon\leq\varepsilon_0.
	\end{equation*}
Consequently, for fixed $N$, $f_N$, and $\ell$,
$\varepsilon^{-1}\mathcal B_{\ell,\varepsilon}(f_N)
\to\mathcal F_\ell(f_N)$ as $\varepsilon\to0$.
\end{theorem}

\begin{proof}
	Taking $\varphi=L_\ell$ in the kernel expansion
	\eqref{eq:fp_kernel_expansion} gives
	\begin{equation*}
		\frac{1}{\varepsilon}
		\mathbb{E}_{\mu}\!\left[
		L_\ell(w_\varepsilon')-L_\ell(w)
		\right]
		=
		b(w)(w_{\ast}-w)L_\ell'(w)
		+
		\frac{\lambda}{2}d(w)^2L_\ell''(w)
		+
		R_{\ell,\varepsilon}(w,w_{\ast}),
	\end{equation*}
	where $R_{\ell,\varepsilon}:=R_{L_\ell,\varepsilon}$
	satisfies
	$\|R_{\ell,\varepsilon}\|_{L^\infty(I^2)}
	\leq C_\ell\sqrt{\varepsilon}$
	by the estimate in the proof of
	Proposition~\ref{prop:quasi_invariant_limit}.
	
	Using $\int_I f_N\,\mathrm dw=1$ and
	$\int_I wf_N\,\mathrm dw=m_N$, we obtain
	\begin{equation*}
	\frac{1}{\varepsilon}\mathcal B_{\ell,\varepsilon}(f_N)
	-\mathcal F_\ell(f_N)
		=
		\int_{I^2}
		R_{\ell,\varepsilon}(w,w_{\ast})
		f_N(w)f_N(w_{\ast})
		\,\mathrm dw\,\mathrm dw_{\ast}.
	\end{equation*}
	Therefore,
	\begin{equation*}
		\begin{aligned}
			\left|
			\frac{1}{\varepsilon}\mathcal B_{\ell,\varepsilon}(f_N)
			-\mathcal F_\ell(f_N)
			\right|
			&\leq
			\|R_{\ell,\varepsilon}\|_{L^\infty(I^2)}
			\int_{I^2}
			|f_N(w)f_N(w_{\ast})|
			\,\mathrm dw\,\mathrm dw_{\ast}\\
			&\leq
			C_\ell\|f_N\|_{L^1(I)}^2\sqrt{\varepsilon}.
		\end{aligned}
	\end{equation*}
	This proves the error estimate, and the convergence follows by
	letting $\varepsilon\to0$.
\end{proof}

Therefore, for fixed $N$, the quasi-invariant limit and the
Legendre--Galerkin discretization commute at the level of the exact Galerkin operators.

\section{Numerical experiments}
\label{sec:numerical_experiments}

In the following, several numerical experiments are presented to support the theoretical results. After that, we try to simulate
the formation of opinion clusters under the bounded-confidence interactions.

\subsection{Verification of theoretical results}
\label{subsec:accuracy_conservation}
We use
$\alpha(w)=\frac{3}{4}+\frac{1}{4}w^2$,
$\beta(w)=\frac{1}{4}- \frac{1}{4}w^2$,
$d(w)=1-w^2$,
$\eta=0.05U$,
$U\sim\mathcal U[-1,1]$, \(\nu=1\), with initial density
\begin{equation}
	f_0(w)
	=
	\frac{1}{Z_0}\exp\bigl(-10(w+0.25)^2\bigr),
	\qquad
	Z_0
	=
	\int_I\exp\bigl(-10(w+0.25)^2\bigr)\,\mathrm dw.
	\label{eq:numerical_smooth_initial_datum}
\end{equation}
For this interaction,
$\partial_w w'\geq0.63$ uniformly for
$w,w_\ast\in I$ and $|\eta|\leq0.05$.
Thus, \eqref{eq:error_collision_bilinear} follows
by changing variables in $w$ and applying the
same duality argument as above.
For a numerical solution \(f_{\rm num}\) and a reference solution
\(f_{\rm ref}\) at the final time, we define the relative \(L^2(I)\)-error
by
$E_{\rm rel}
:=
{\|f_{\rm num}-f_{\rm ref}\|_{L^2(I)}}/
{\|f_{\rm ref}\|_{L^2(I)}}$.

For the spatial convergence test, we take \(T=0.2\) and
\(\Delta t=10^{-4}\). The reference solution is computed with
\(N_{\rm ref}=128\) and \(N_{q,\rm ref}=264\), whereas the tested
solutions use \(N_q=2N+8\). The results are reported in Table~\ref{tab:spatial_convergence}.
For the temporal convergence test, we fix \(N=32\), \(N_q=72\), and
\(T=2\). The reference solution is computed with
\(\Delta t_{\rm ref}=2^{-11}\), and the observed order is defined by
$\log_2\left(
\frac{E_{\rm rel}(2\Delta t)}
{E_{\rm rel}(\Delta t)}
\right)$. The results are reported in
Table~\ref{tab:temporal_convergence}.
Mass conservation is also verified numerically.

\begin{table}[h!]
	\centering
	\caption{Spatial convergence.}
	\vspace{-5pt}
	\label{tab:spatial_convergence}
	\renewcommand{\arraystretch}{1.2}
	\begin{tabular}{c|c|c|c}
		\hline\hline
		$N$ & $N_q$ &$E_{\rm rel}$ & Convergence  \\ \hline
		8  & 24& $ 8.199\times 10^{-2}$ & \multirow{5}{*}{\includegraphics[width=0.26\textwidth]{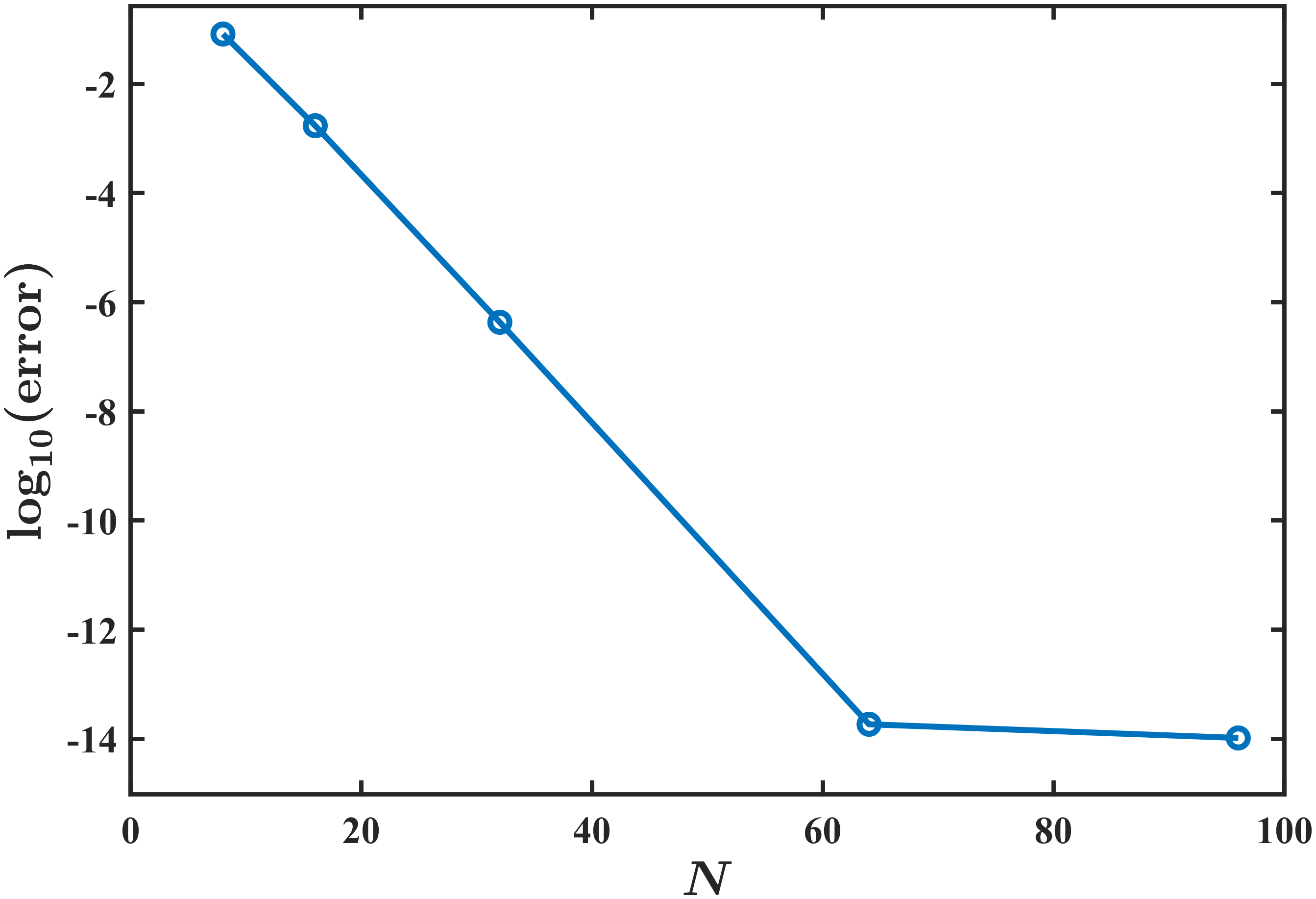}} \\
		16 & 40 &$ 1.713 \times 10^{-3}$ & \\
		32 &  72&$ 4.286 \times 10^{-7}$ & \\
		64 & 136&$ 1.845 \times 10^{-14}$ & \\
		96 & 200&$ 1.037 \times 10^{-14}$ & \\ \hline\hline
	\end{tabular}
\end{table}

\begin{table}[h!]
	\centering
	\caption{Temporal convergence.}
	\vspace{-5pt}
	\label{tab:temporal_convergence}
	\renewcommand{\arraystretch}{1.25}
	\begin{tabular}{c|c|c|c}
				\hline\hline
				$\Delta t$ & $E_{\rm rel}$ & Order & Convergence  \\ \hline
				$0.25$ & $5.207 \times 10^{-7}$ & - & \multirow{5.2}{*}{\includegraphics[width=0.26\textwidth]{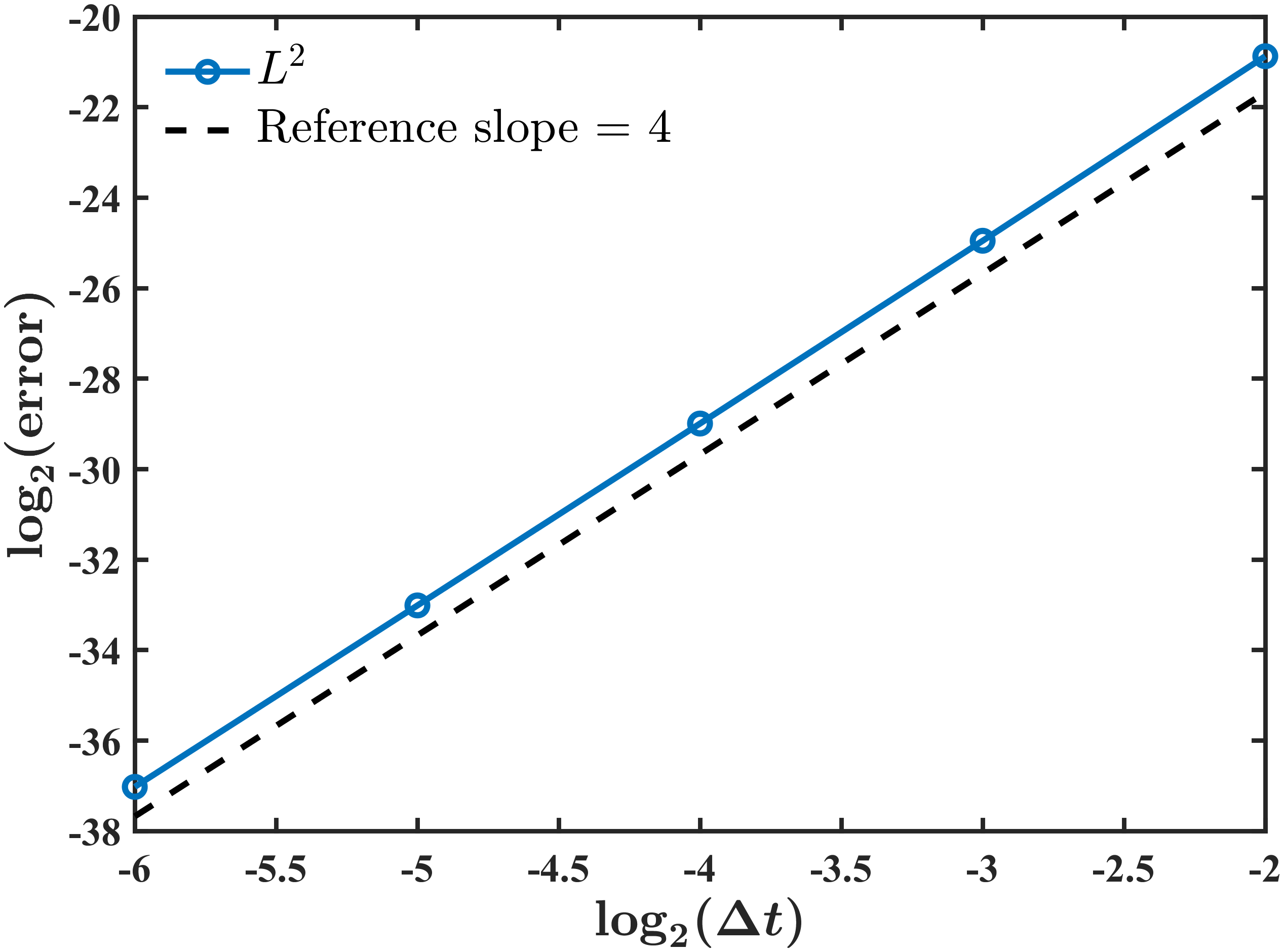}} \\
				$0.125$ & $3.081 \times 10^{-8}$ & 4.079& \\
				$0.0625$ & $1.873 \times 10^{-9}$ & 4.040 & \\
				$0.03125$ & $1.155 \times 10^{-10}$ & 4.020 & \\ 
				$0.015625$ & $7.169 \times 10^{-12}$ & 4.010 & \\
				\hline\hline
			\end{tabular}
\end{table}

\subsection{Discrete quasi-invariant consistency}
\label{subsec:qi_consistency}

We numerically verify the discrete quasi-invariant consistency
established in Theorem~\ref{thm:qi_discrete_consistency} for the
scaled interaction \eqref{eq:qi_scaled_interaction}. Here we choose
$b(w)=\frac14(1-w^2)$, $d(w)=1-w^2$,
$\lambda=5\times10^{-2}$, and $\nu=1$. The random variable $\eta$ is
defined by $\mathbb P(\eta=1/\sqrt2)=2/3$ and
$\mathbb P(\eta=-\sqrt2)=1/3$. Hence,
$\mathbb E_{\mu}[\eta]=0$,
$\mathbb E_{\mu}[\eta^2]=1$, and
$\mathbb E_{\mu}[\eta^3]=-1/\sqrt{2}$.

\begin{figure}[h!]
	\centering
	\includegraphics[width=0.4\linewidth]
	{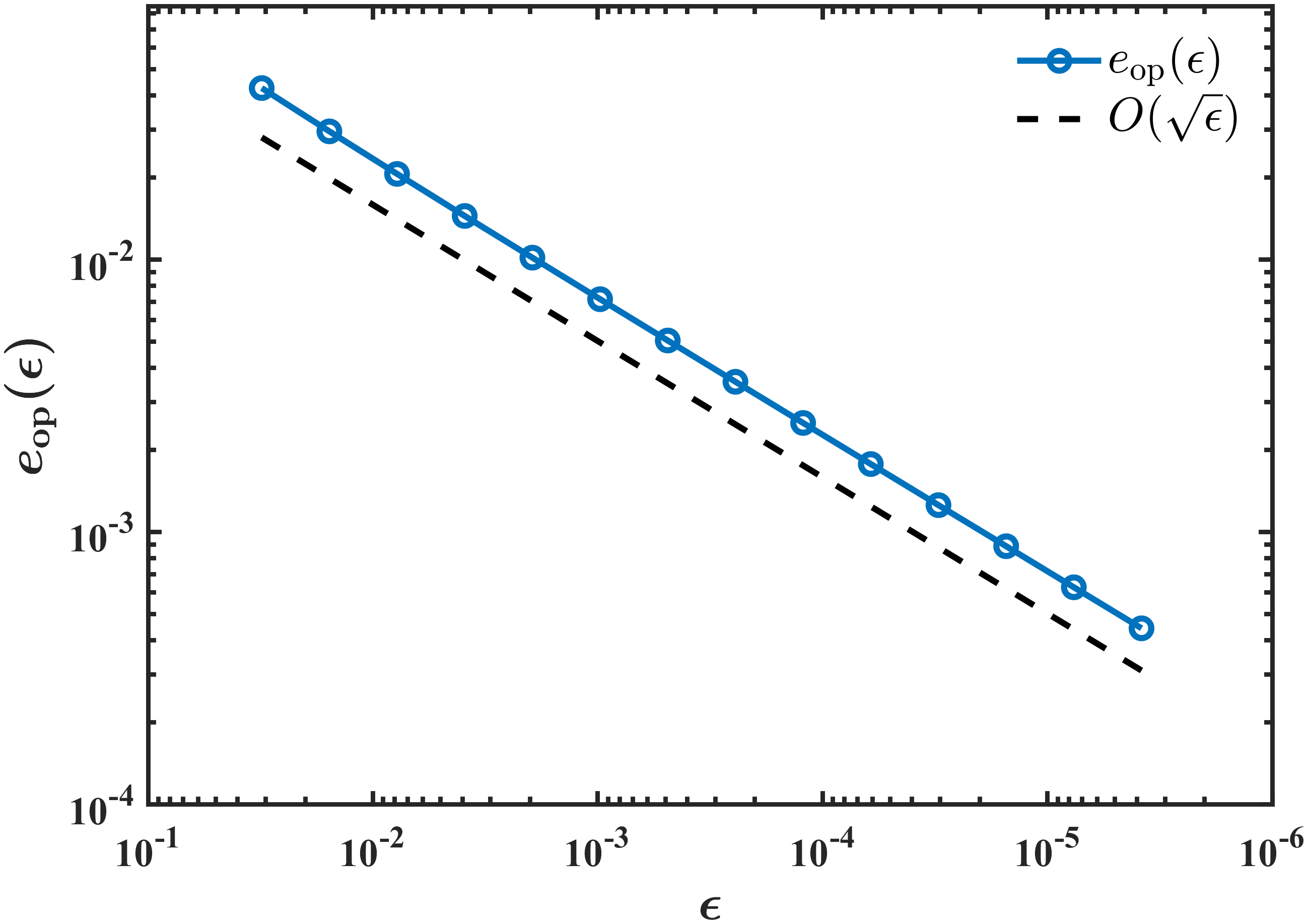}
	\caption{Discrete quasi-invariant consistency is verified when $\epsilon$ is sufficiently small.}
	\label{fig:consistency}
\end{figure}

We use the smooth initial density
\eqref{eq:numerical_smooth_initial_datum} and consider
$0<\varepsilon\leq2^{-5}$, for which the admissibility condition
\eqref{eq:qi_admissibility} is satisfied. Let $f_N^0=\Pi_N f_0$ denote the Legendre projection
of the initial density $f_0$. We define the discrete
operator error by
$e_{\mathrm{op}}(\varepsilon)
:=
\left(
\sum_{\ell=0}^{N}
\frac{1}{h_\ell}
\left|
\frac{1}{\varepsilon}
\mathcal B_{\ell,\varepsilon}(f_N^0)
-
\mathcal F_\ell(f_N^0)
\right|^2
\right)^{1/2}$. 
Figure~\ref{fig:consistency} verifies that \(e_{\rm op}(\varepsilon)=\mathcal O(\sqrt{\varepsilon})\) under $N=24$ and
$N_q=64$, which coincides with Theorem~\ref{thm:qi_discrete_consistency}.

\subsection{Bounded-confidence clustering}
\label{subsec:bounded_confidence}

Now consider \eqref{eq:bounded_confidence_interaction}
with $\theta(w)\equiv\gamma$ and apply the scheme
described in Section~\ref{subsec:bounded_confidence_method}.
We take \(f_0(w)=1/2\), \(d(w)=1-w^2\), and use 
the truncated Gaussian noise
\(\eta=\sigma Z/\sqrt{v_c}\), where
\(Z\sim\mathcal N(0,1)\) conditioned on \(|Z|\leq c\) and
\(v_c=1-2c\varphi(c)/(2\Phi(c)-1)\), with \(\varphi\) and \(\Phi\)
denoting the standard normal density and cumulative distribution
function, respectively.
This normalization gives
$\mathbb E_\mu[\eta]=0$ and
$\mathbb E_\mu[\eta^2]=\sigma^2$, while
$\eta_{\max}=\sigma c/\sqrt{v_c}$.

\begin{figure}[h!]
	\centering
	\includegraphics[width=0.32\linewidth]
	{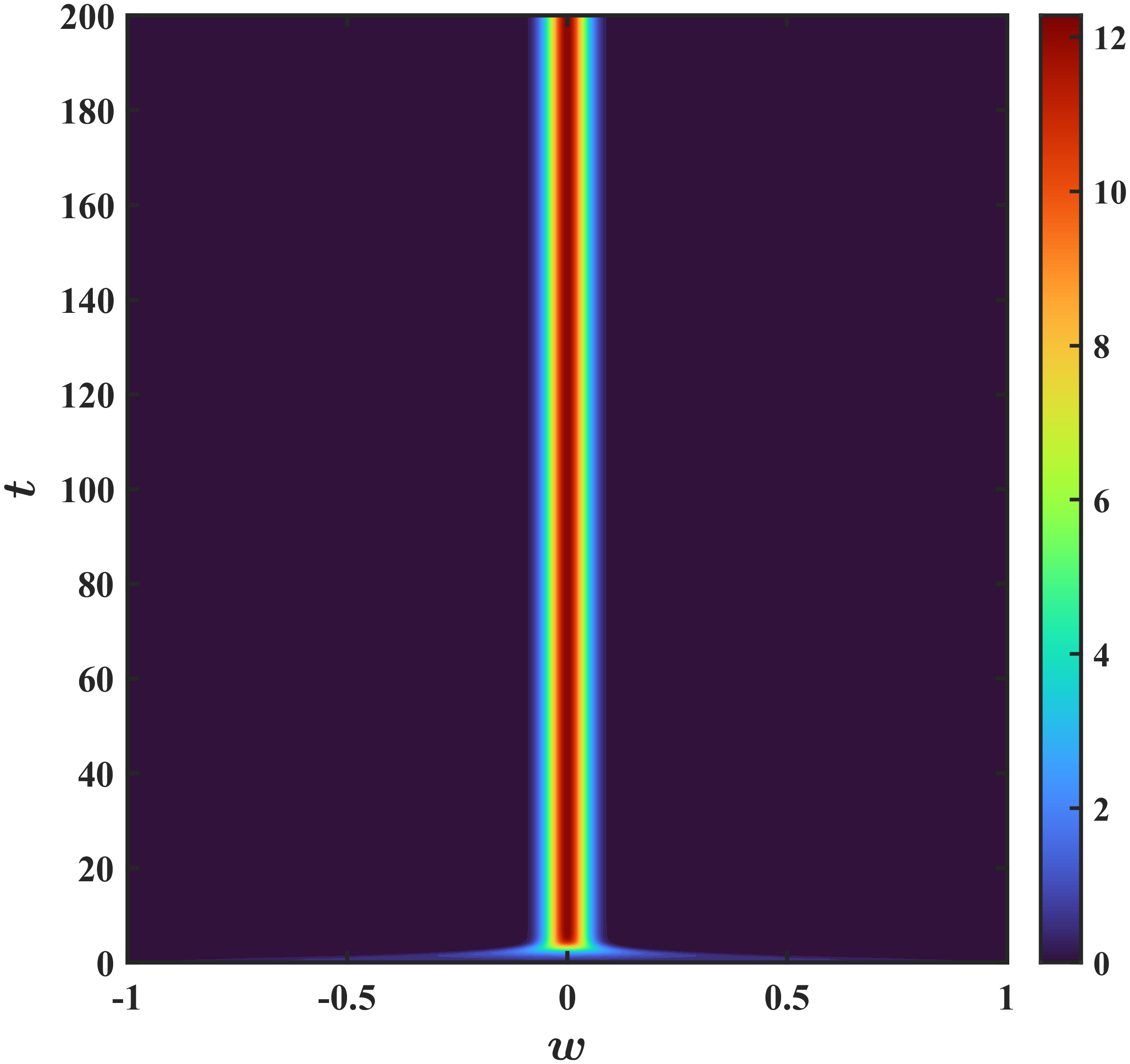}
	\includegraphics[width=0.32\linewidth]
	{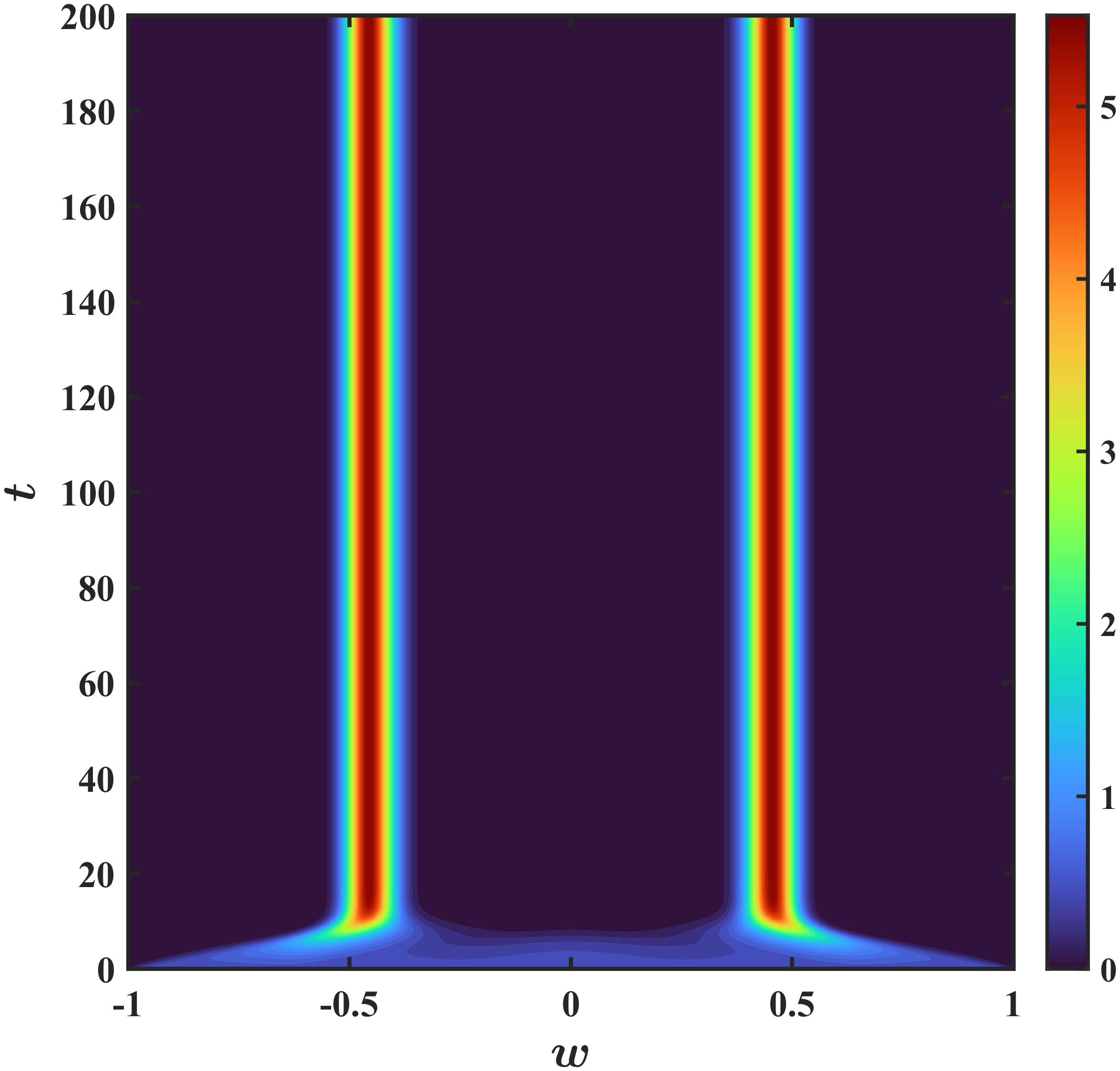}
	\includegraphics[width=0.32\linewidth]
	{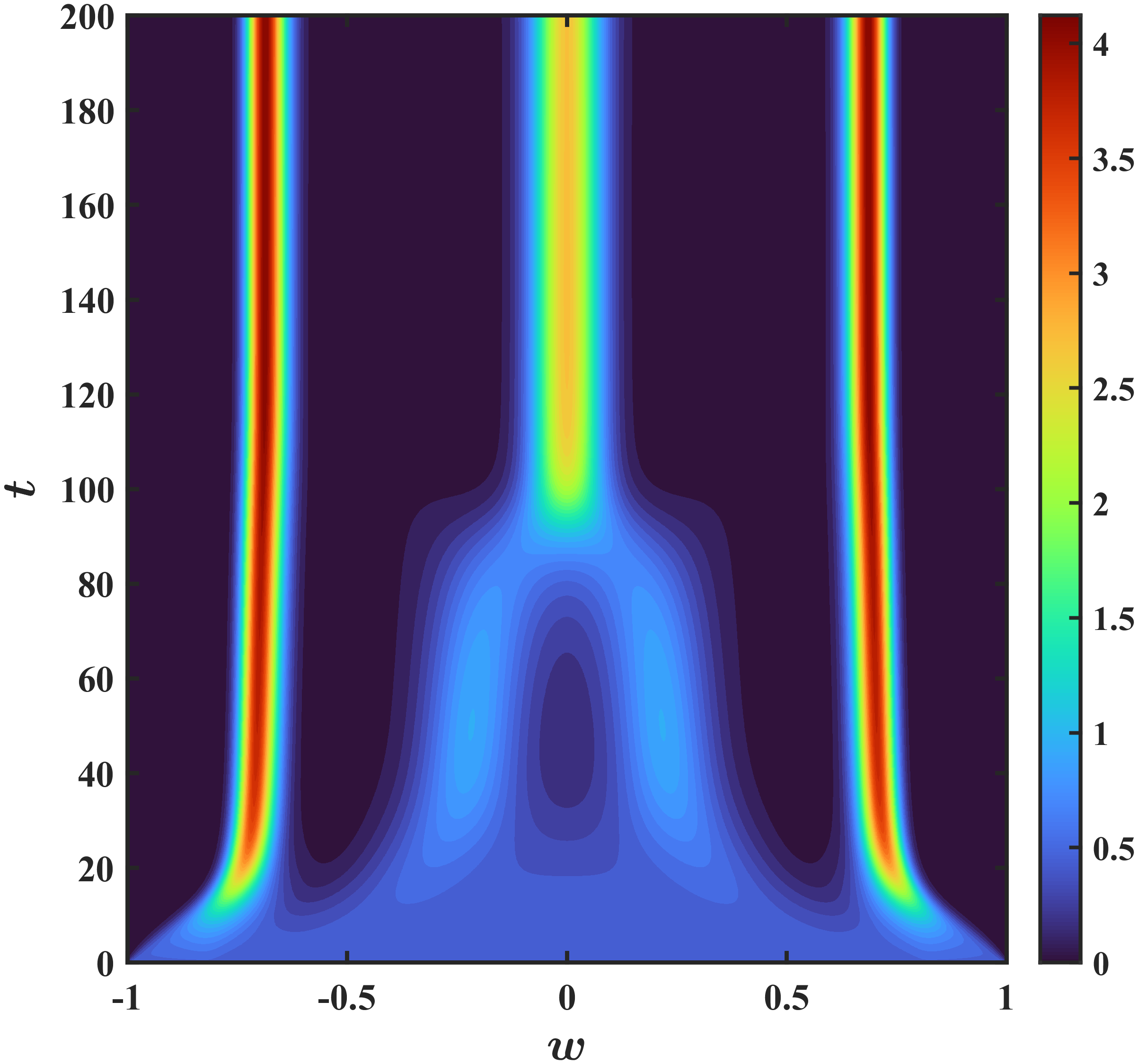}
	\includegraphics[width=0.32\linewidth]
	{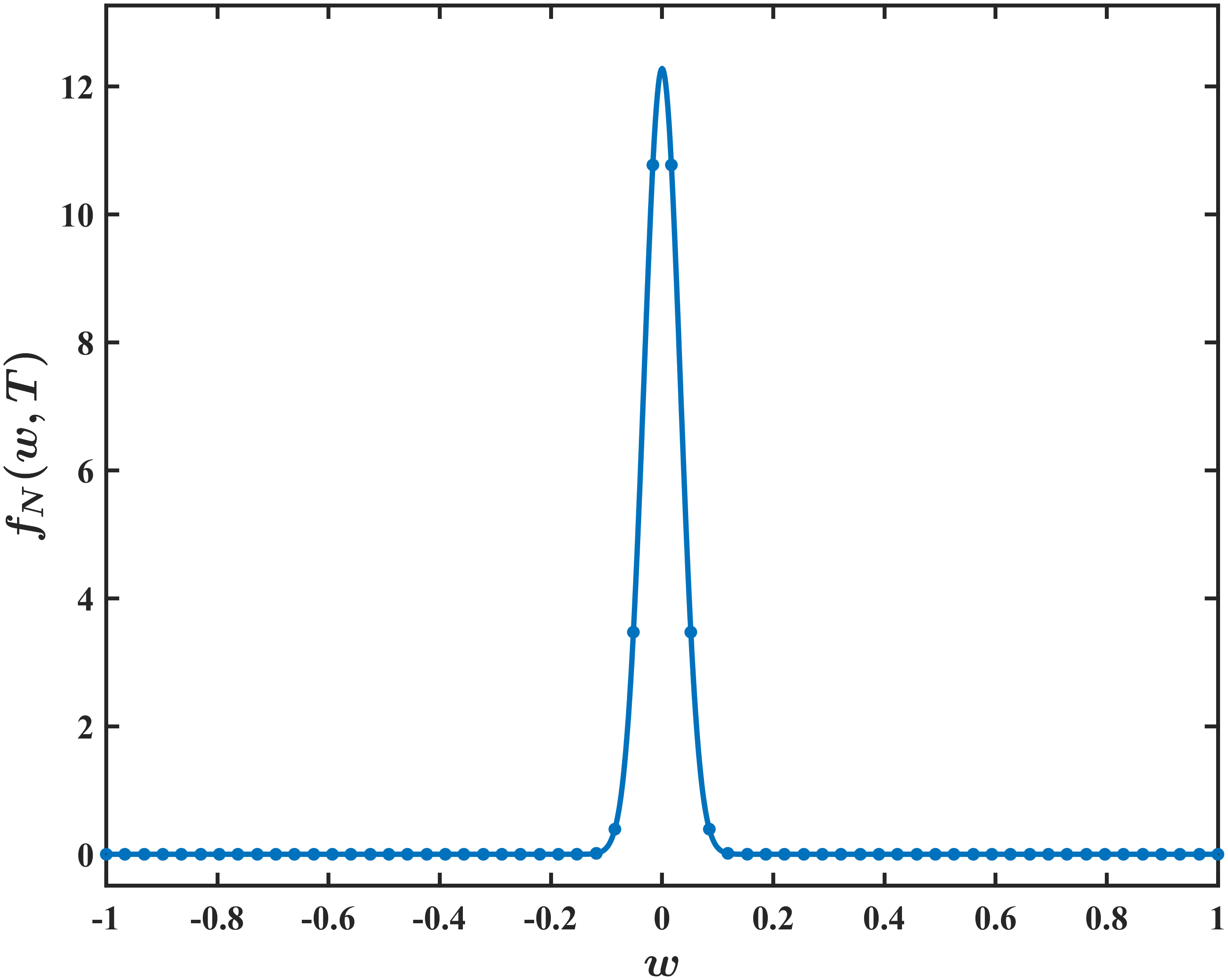}	
	\includegraphics[width=0.32\linewidth]
	{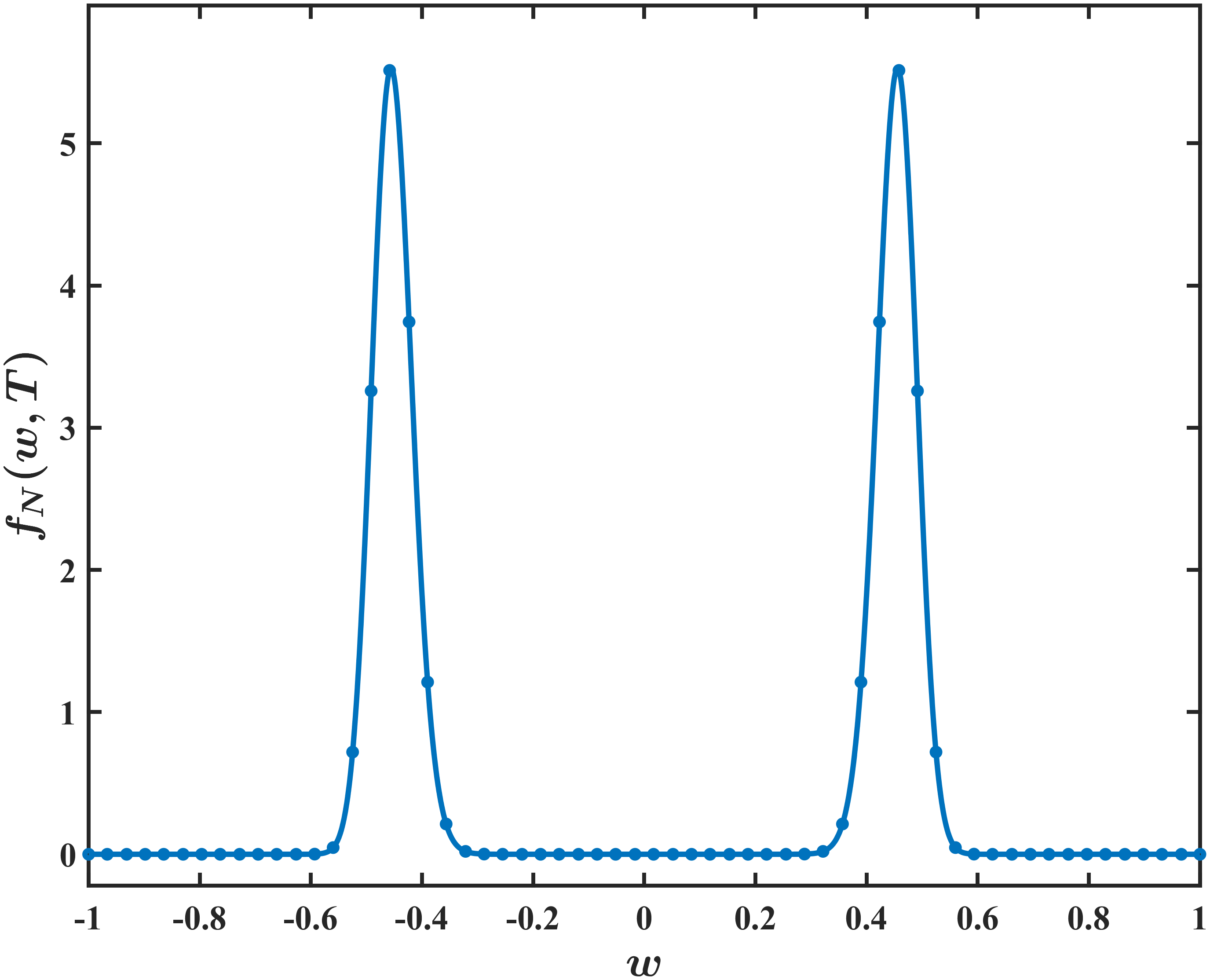}	
	\includegraphics[width=0.32\linewidth]
	{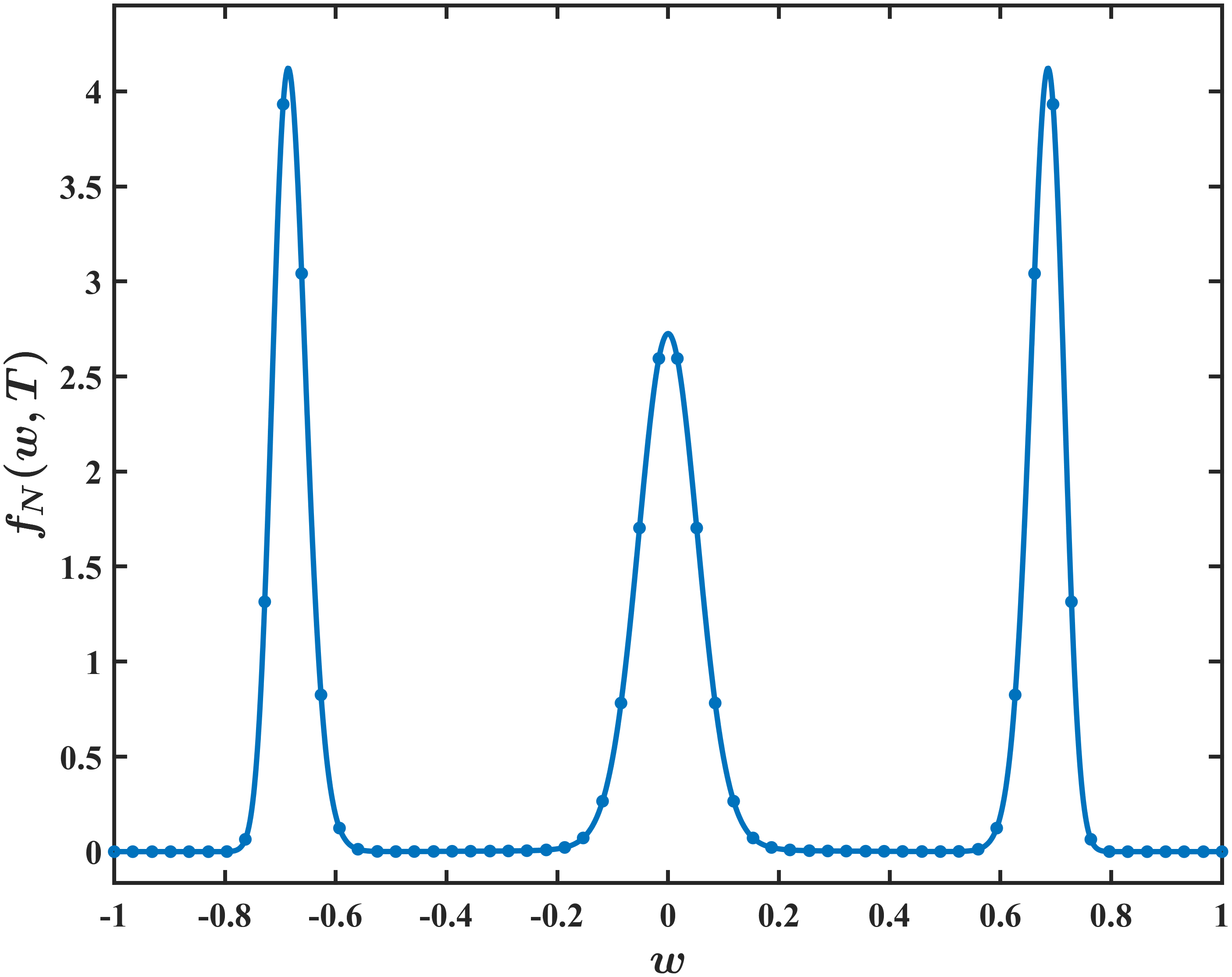}
	\caption{Bounded-confidence clustering for different confidence
		thresholds (from left to right: \(\Delta=1\), \(0.4\), and \(0.2\)).
		The top row shows the evolution of the opinion density, while the bottom row shows the corresponding
		final-time profiles \(f_N(w,T)\). More clusters are formed for smaller threshold $\Delta$, indicating that the polarization of opinions may occur when people tend to interact with smaller neighbors.}
	\label{fig:bounded_confidence_clustering}
\end{figure}

The parameters are set as \(\gamma=5\times10^{-2}\),
\(\sigma=10^{-2}\), \(c=3\), \(\nu=20\), \(N=160\),
\(\Delta t=2.5\times10^{-3}\), and \(T=200\). We use $N_q=200$ Gauss--Legendre nodes on each
nondegenerate subinterval. The required noise
moments are evaluated analytically, and the sufficient admissibility
condition \(2\eta_{\max}\leq1-\gamma\) is satisfied. We consider
\(\Delta=1\), \(0.4\), and \(0.2\).
Figure~\ref{fig:bounded_confidence_clustering} shows that decreasing
\(\Delta\) localizes the interactions and leads to the formation of
separated opinion clusters, which is also demonstrated by movies in our supplementary materials.
This indicates that the polarization of opinions may occur when people tend to interact with smaller neighbors.

\section{Conclusion}
\label{sec:conclusion}
In this paper, we proposed an efficient Legendre-Galerkin spectral method for a class
of Boltzmann-type opinion equations with opinion-dependent
interaction rules.
An exact finite decoupling expansion of the
weak collision kernel leads to a tensor representation of the
Galerkin collision operator.
A recursive construction using noise moments and the
Legendre three-term recurrence allows the collision tensor to be precomputed via
1D Gauss-Legendre quadrature, thereby reducing each collision evaluation to tensor contractions.
The scheme can be extended to bounded-confidence interactions.

Under the stated regularity and stability assumptions,
we established an $L^2$-error estimate.
In the quasi-invariant regime, we derived the corresponding
Fokker--Planck equation and established consistency between
the exact Legendre--Galerkin operators of the rescaled
Boltzmann equation and its limit for each fixed polynomial
degree. Numerical experiments supported the theoretical accuracy, mass-conservation, and the quasi-invariant
consistency. Moreover, the scheme can accurately capture the formation of opinion clusters under the bounded-confidence interactions.

The proposed spectral solver for opinion dynamics also provides a solid foundation for the model parameter inferences and may facilitate a data-driven approach to studying the collective opinion dynamics, which will be discussed in our future work.


\section*{Acknowledgments}
The authors assume responsibility for all content.

\bibliographystyle{siamplain}

\end{document}